\documentclass [12pt]{article}
\usepackage{amsfonts,amsthm,amsmath,amssymb}
\usepackage{amscd,graphicx,psfrag,subfigure}
\usepackage{epsfig}

\newcommand{\bn}{{\mathbb N}}
\newcommand{\bc}{{\mathbb C}}
\newcommand{\br}{{\mathbb R}}
\newcommand{\bz}{{\mathbb Z}}
\newcommand{\bbZ}{{\Bbb Z}}

\newcommand{\bbR}{{\Bbb R}}
\newcommand{\bbC}{{\Bbb C}}
\newcommand{\bbN}{{\Bbb N}}
\newcommand{\bbD}{{\Bbb D}}
\newcommand{\T}{{\mathfrak T}}

\newcommand{\J}{{\mathcal J}}
\newcommand{\F}{{\mathcal F}}
\newcommand{\A}{{\mathcal A}}

\usepackage{color}
\definecolor{green}{rgb}{0,0.7,0.3}
\definecolor{red}{rgb}{0.9,0,0.5}

\newtheorem{Theorem}{Theorem}[section]
\newtheorem{Corollary}[Theorem]{Corollary}
\newtheorem{Definition}[Theorem]{Definition}
\newtheorem{Lemma}[Theorem]{Lemma}
\newtheorem{Proposition}[Theorem]{Proposition}
\newtheorem{Remark}[Theorem]{Remark}

\newsavebox{\savepar}

\title{Non-landing hairs in Sierpi\'nski curve\\ Julia sets
of transcendental entire maps\\ (Revised Version)}

\author{{\small Antonio Garijo} \and {\small Xavier
Jarque} \and {\small M\'{o}nica
Moreno Rocha\footnote{Corresponding author. E-mail:~\texttt{mmoreno@cimat.mx}.}}}

\date{\today}

\begin{document}
\maketitle

\abstract{ We consider the family of transcendental entire maps given by $f_a(z)=a(z-(1-a))\exp(z+a)$ where $a$ is a complex parameter. Every map has a superattracting fixed point at $z=-a$ and an asymptotic value at $z=0$. For $a>1$ the Julia set of $f_a$ is known to be homeomorphic to the Sierpi\'nski universal curve~\cite{Moro}, thus containing embedded copies of any one-dimensional plane continuum. In this paper we study subcontinua of the Julia set that can be defined in a combinatorial manner.
In particular, we show the existence of non-landing hairs with prescribed combinatorics  embedded in the Julia set for all parameters $a\geq 3$.  We also study the relation between non-landing hairs and the immediate basin of attraction of $z=-a$. Even as each non-landing hair accumulates onto the boundary of the immediate basin at a single point, its closure, nonetheless,  becomes an indecomposable subcontinuum of the Julia set.}

\vfill
\bigskip
\noindent
{\bf Keywords:} Transcendental entire maps, Julia set, non-landing hairs, indecomposable continua.

\bigskip
\noindent
{\bf Mathematics Subject Classification (2000):} 37F10, 37F20.

\pagebreak

\section{Introduction}\label{section:intro}

Let $f:\bbC \to \bbC$ be a transcendental entire map. The \emph{Fatou set} ${\cal F}(f)$ is the largest open set where iterates of $f$ form a normal family. Its complement in $\bbC$ is the \emph{Julia set} ${\cal J}(f)$ and it is a non-empty and unbounded subset of the plane. When the set of singular values  is  bounded, we say $f$ is of \emph{bounded singular type} and denote this class of maps by $\cal B$. It has been shown in \cite{Ba} and \cite{R} that the Julia set of a hyperbolic map in $\cal B$ contains uncountably many unbounded curves, usually known as \emph{hairs}, \cite{DT}.~A hair is said to \emph{land} if it is homeomorphic to the half-closed ray $[0,+\infty)$. The point corresponding to $t=0$ is known as the \emph{endpoint} of the hair. In contrast, if its accumulation set is a non-trivial continuum, we obtain a \emph{non-landing} hair.

In this paper we study a particular class of non-landing hairs in the Julia set of transcendental entire maps given by
\[\label{eq:themap}
f_a(z)=a (z-(1-a)) \exp(z+a),
\]
when $a$ is a real parameter. For all complex values of $a$, the map $f_a$ has a superattracting fixed point at $z=-a$ and an asymptotic value at the origin whose dynamics depends on the parameter $a$. If the orbit of the asymptotic value escapes to $+\infty$, we say $a$ is an {\it escaping parameter}. For example, when $a>1$, the orbit of the asymptotic value escapes to $+\infty$ along the positive real axis.

To our knowledge, the family $f_{a}$ was first introduced by Morosawa, \cite{Moro}, as an example of a transcendental entire map whose Julia set is homeomorphic to the {\it Sierpi\'nski curve} continuum when $a>1$.  Any planar set that is compact, connected, locally connected, nowhere dense, and has the property that any two complementary domains are bounded by disjoint simple closed curves is homeomorphic to the Sierpi\'nski curve continuum (Whyburn, \cite{Why}). It is also a \emph{universal} continuum, in the sense that it contains a homeomorphic copy of every one-dimensional plane continuum (Kuratowski, \cite{K}). We take advantage of this property to combinatorially construct subsets of $\J(f_a), \ a >1$,  that in turn are  \emph{indecomposable  continua}. An {\it indecomposable continuum} is a compact, connected set that cannot be written as the union of two proper connected and closed subsets.~Observe that a landing hair together with the point at infinity is in fact a decomposable continuum.

Every known example in the literature of indecomposable subsets of Julia sets arises from a single family of maps, namely the exponential family $E_{\lambda}(z)=\lambda \exp(z)$.  The first example was given by Devaney \cite{D} when $\lambda=1$ so the asymptotic value escapes to infinity and the Julia set is the whole plane. Under the assumption that either the asymptotic value escapes to infinity or has a preperiodic orbit (thus holding again ${\cal J}(E_\lambda)= \bbC$), several authors have been able to construct topologically distinct indecomposable continua embedded in $\J(E_\lambda)$ (see among other works, \cite{DJ1}, \cite{DJM}, and \cite{R1}, where a generalization of previous results for a large set of $\lambda$-parameters can be found). 

\begin{figure}[hbt]
    \centering
    \subfigure[\scriptsize{Parameter~$a = 1.1$. } ]{
    \includegraphics[width=0.3 \textwidth]{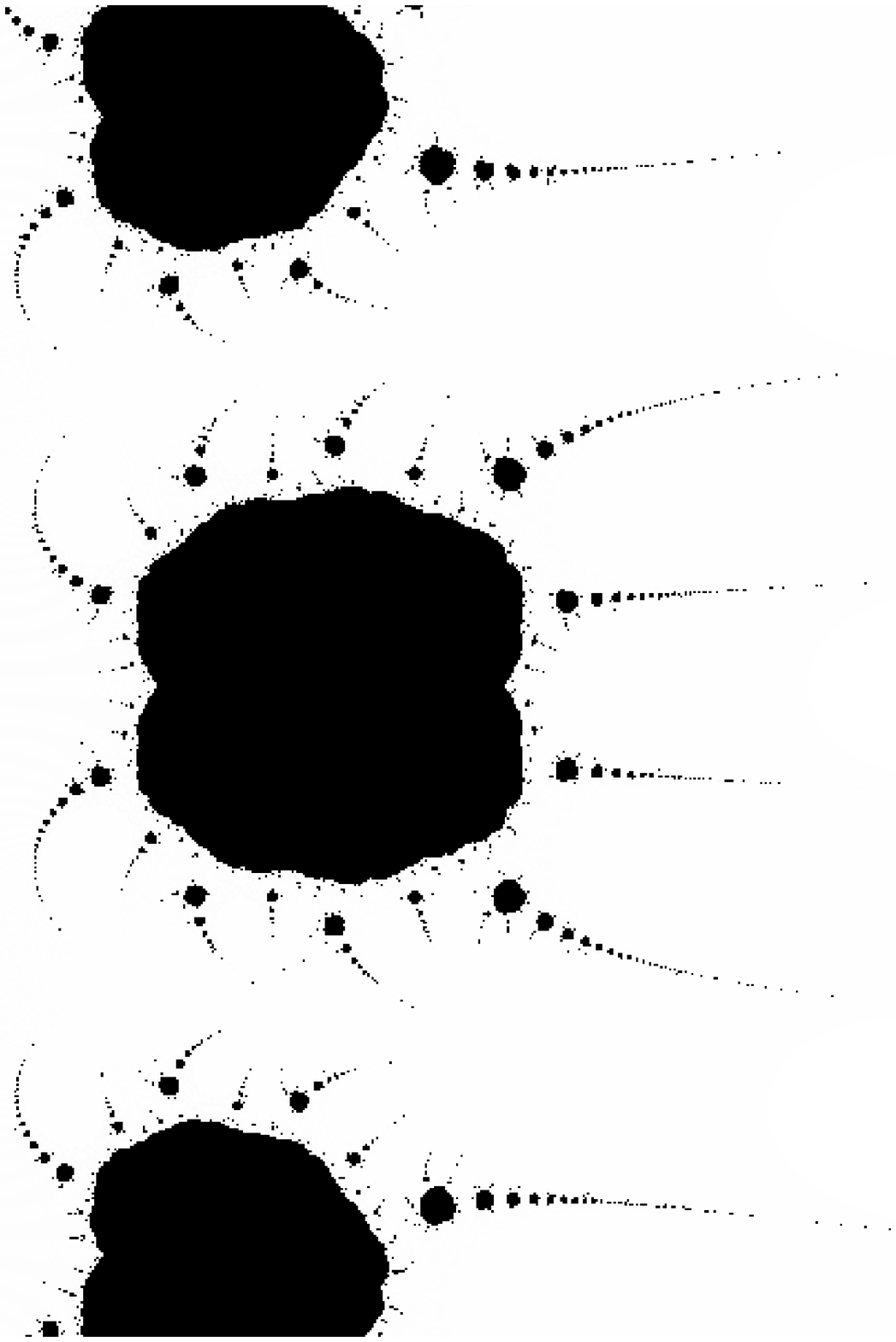}}
    \hspace{0.0in}
    \subfigure[\scriptsize{Parameter  $a =2.1$.}  ]{
    \includegraphics[width=0.3 \textwidth]{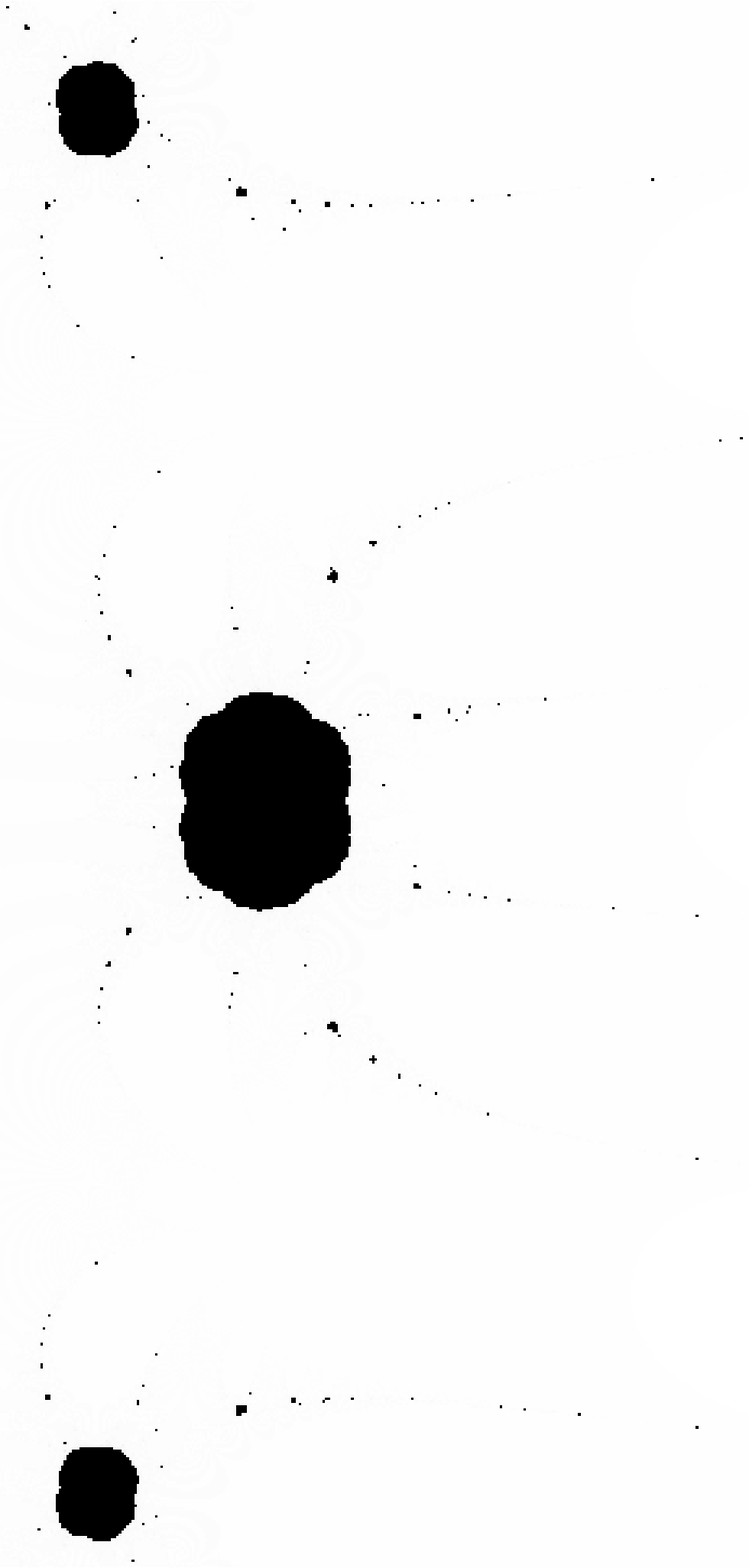}}
    \hspace{0.0in}
    \subfigure[\scriptsize{Parameter  $a =3.1$. }  ]{
     \includegraphics[width=0.3 \textwidth]{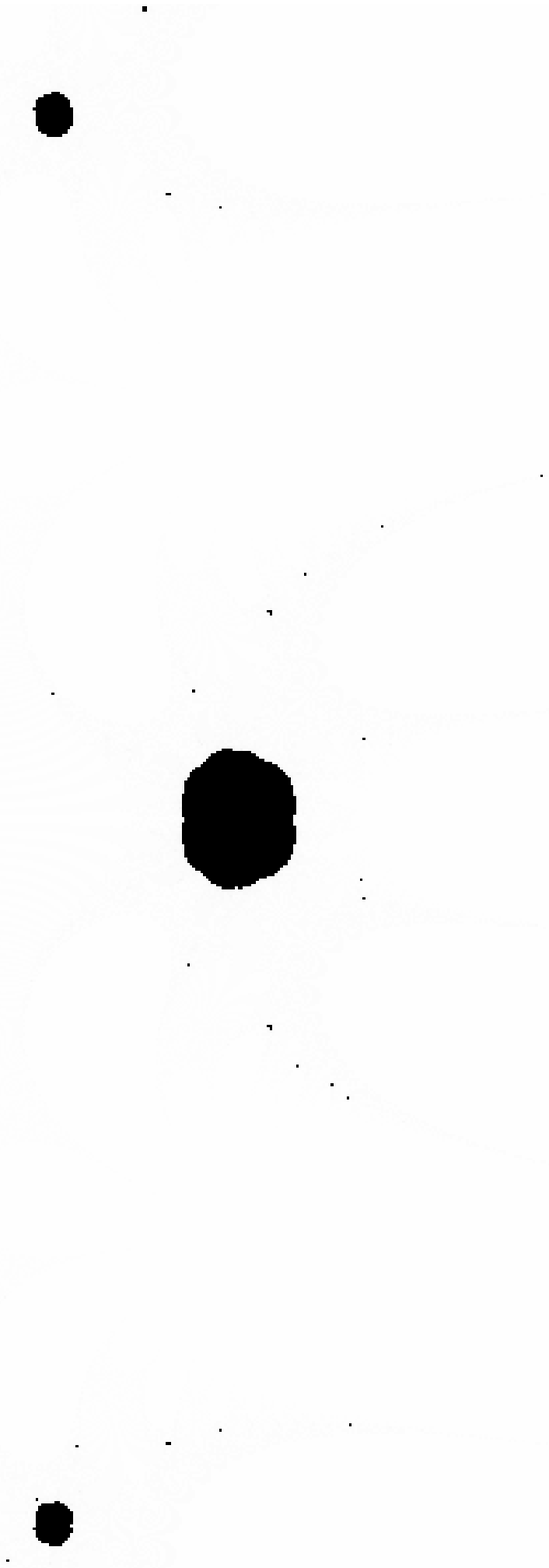}}
    \caption{\small{The Julia set for $f_a$ (and $a$ an escaping parameter) is shown in white, the Fatou set, in black. }}
    \label{fig:julia_set}
\end{figure}

Our work provides examples of indecomposable subcontinua of Julia sets outside the exponential family and without the assumption that $\J(f_a)$ equals $\bbC$, since $f_a$ has a superattracting fixed point  for all $a\in\mathbb C$.
Denote by $\A(-a)$ the \emph{basin of attraction} of $-a$, that is, the set of points with forward orbits converging to $-a$. Also denote by $\mathcal{A}^*(-a)$ the \emph{immediate basin of attraction} of $-a$ which is  the connected component of $\mathcal{A}(-a)$ containing $-a$. In \cite{Moro} Morosawa showed that all connected components of $\A(-a)$ are bounded Jordan domains.  Moreover, whenever $a>1$, the orbit of the free asymptotic value escapes to infinity. Since there are no other singular values, $\F(f_a)$ cannot contain another attracting basin, or a parabolic basin, or a Siegel disk as these components must be associated with a non-escaping singular value.~Maps with a finite number of singular values  do not exhibit neither wandering domains (\cite{EL2,GK})   nor Baker domains (\cite{EL2}). Hence $\mathcal{F}(f_{a})=\mathcal{A}(-a)$. In Figure \ref{fig:julia_set}, we display the dynamical plane of $f_a$ for different values of $a>1$. The basin of attraction of $-a$ is shown in black, while points in the Julia set are shown in white. 

Let us summarize our main results. Since $\J(f_a)$ is homeomorphic to the Sierpi\'nski universal curve, it must contain embedded copies of planar indecomposable continua, so we obtain some of them in terms of its combinatorics. To do so, we first characterize the topology and dynamics of the boundary of $\mathcal{A}^*(-a)$ by a polynomial-like construction (Proposition \ref{proposition:pol_like}). Then, using general results of transcendental entire maps, we obtain curves in the Julia set contained in the far right plane and with specific combinatorics (Proposition~\ref{prop:ConjugTails}). By a controlled process of consecutive pullbacks of some of these curves, we extend them into non-landing hairs that limit upon themselves at every point (Theorem \ref{theorem:indecom}). Using a result due to Curry, \cite{C}, we show the closure of such hairs are indecomposable continua (Theorem \ref{thm:DoesNotSeparate}). Finally, we study the relation between each indecomposable continuum and the boundary of $\A^*(-a)$ showing that the intersection between these two sets reduces to a unique point (Theorem \ref{theorem:relation_inde_basin}). 
As a consequence of these results, we show the existence of a dense set of points in $\partial \mathcal{A}^*(-a)$ that are landing points of a unique hair (in particular there are no \emph{pinchings} that arise as in other maps in class $\cal B$ having a superattracting basin), while there is a residual set of points in $\partial \mathcal{A}^*(-a)$ that, even though they belong to the accumulation set of a certain ray, they are not landing points of hairs.

The outline of this paper is as follows: in \S
\ref{section:dyn_plane} we describe the dynamical plane of $f_a$ for $a\geq 3$.  \S \ref{section:targets} contains most of our technical results while in \S \ref{section:indecom} we provide the proofs of our main results.

\begin{figure}[hbt]
  
  \psfrag{a}[][]{\footnotesize $\zeta_1$} 
  \psfrag{b}[][]{\footnotesize $\eta_1$} 
  \psfrag{d}[][]{\footnotesize $\zeta_{-1}$} 
  \psfrag{c}[][]{\footnotesize $\eta_{-1}$} 
  \psfrag{e}[][]{\footnotesize $T_1$} 
  \psfrag{f}[][]{\footnotesize $T_0$} 
  \psfrag{r}[][]{\footnotesize $R$} 
  \psfrag{h}[][]{\footnotesize $H_R$} 
  \centerline{\includegraphics[width=0.5\textwidth]{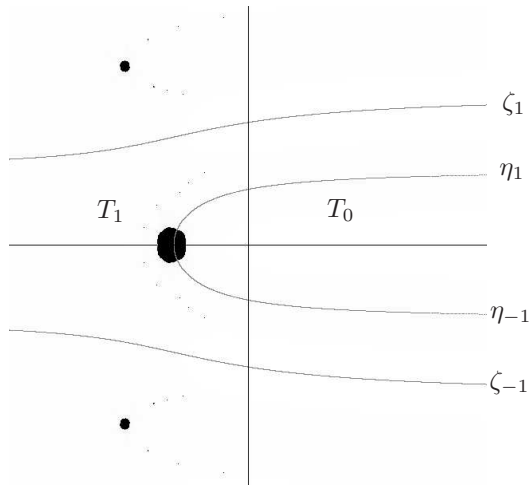}}
   \caption{\small{Dynamical plane of $f_{3.1}$ and regions $T_0\cup T_1$. The boundary curves $\zeta_1$ and $\zeta_{-1}$ extend from $-\infty$ to $+ \infty$ while curves $\eta_1$ and $\eta_{-1}$ meet at $-a=-3.1$ and extend to $+\infty$.}}
      \label{fig:sigmas}
\end{figure}

\noindent
\emph{Notation and terminology.}

$B_\varepsilon(x)=\{z\in \bbC~|~|z-x|<\varepsilon\}$.

$\overline{U}$ denotes the closure of a set $U$.

Connected components will be refered to as components.

A curve $\gamma$ \emph{cuts across} a
\begin{enumerate}
\item line $L$ if the intersection $\gamma\cap L$ is not tangential,
\item rectangle $R$ if $\gamma$ cuts across both vertical boundaries of $R$ so $\gamma\cap R$ contains a component with endpoints joining those sides,
\item semi-annular region $A$ if $\gamma$ cuts across the inner and outer semicircular boundaries of $A$ so $\gamma\cap A$ contains a connected component with endpoints joining those boundaries. 
\end{enumerate}

\vfill
\pagebreak

\section{Dynamical plane  for escaping real parameters}\label{section:dyn_plane}

Consider escaping parameters of the form $a>1$ for the family of transcendental entire maps
$$f_a(z) = a(z-(1-a))\exp(z+a),$$
which have a unique asymptotic value at $z=0$ and a superattracting fixed point at $z=-a$. For $a>1$, the asymptotic value escapes to infinity along the positive real line and the Fatou set reduces to the basin of attraction of $-a$, $\A(-a)$.  Our first aim in this section is to provide a partition of the complex plane that will allow us to combinatorially analyze the dynamics of points in the Julia set.\\

We start by taking preimages of the forward invariant set $\br^+$. Any point $z=x+iy$ in the complex plane whose image under $f_a$ is a real positive number must satisfy
\begin{equation}\label{eq:def_sigmas}
\begin{split}
& \left( x- (1-a) \right)\cos y  - y\sin y > 0, \\
& \left( x- (1-a) \right)\sin y  + y\cos y = 0. 
\end{split}
\end{equation}

From these conditions, the preimages of  $\br^+$ are infinitely many analytic curves parametrized by $(x,\zeta_j(x))$, with $j\in \bbZ$. For $j=0$, $\zeta_0(x)=0$ and is defined for all  $x \in (1-a,+\infty)$ while the rest of the $\zeta_j$'s are strictly monotonic functions of $x$ defined for all $x \in \br$. When $j\neq 0$, each $\zeta_j$ has two horizontal asymptotes, given by

\[
 \lim\limits_{x\to -\infty}  \zeta_j(x) =   \text{sign}(j) (2|j|-1)\pi i \quad {\rm and} \quad \lim\limits_{x\to +\infty} \zeta_j(x) =2 j \pi  i.
\]

For our purposes, we need to consider the preimage of the interval $(-\infty,-a)$ inside the region bounded by $\zeta_1$ and $\zeta_{-1}$ (see Figure \ref{fig:sigmas}). We obtain two strictly monotonic curves $\eta_1(x)$ and $\eta_{-1}(x)$ defined for $x \in [-a,+\infty)$, satisfying
$$
\lim\limits_{x\to -a^+} \eta_{\pm 1}(x) = 0, \quad  \quad \lim\limits_{x\to +\infty} \eta_{1}(x) = \pi  i  \quad {\rm and} \quad \lim\limits_{x\to +\infty} \eta_{-1}(x) = -\pi i  . 
$$ 

Let $T_{0}$ denote the open and connected set  containing $z=0$ and bounded by  $\eta_1 \cup \eta_{-1}$. Similarly, let $T_{1}$ be the open and connected set bounded by $\zeta_{-1}\cup \eta_{-1} \cup \eta_1 \cup \zeta_{1}$. Far to the right, $T_{1}$ consists of two unbounded and disjoint strips, one above and one below the positive real line.
Since most of our results involve the dynamics of points in $T_0\cup T_1$, we construct a refinement of this region. For $j=0,1$, denote by $T_{j_1}$ and $T_{j_2}$ the proper and disjoint domains in $T_j\setminus \bbR$ with negative and positive imaginary part, respectively.
Finally, for each $j\in \bbZ, j\neq 0,1$, denote by $T_j$ the open and connected strip bounded by the curves $\zeta_{j-1}$ and $\zeta_j$ as Im$(z)$ increases. Then $\{T_j~|~j\in \bbZ \}$ defines the partition of the complex plane sought, while $\{T_{j_i}~|~j=0,1,~i=1,2\}$ defines a refinement of the region $T_0\cup T_1$.\\

It is straightforward to verify that
\begin{equation*}
\begin{split}
& f_a : T_{0} \to \bc \setminus \left(-\infty,-a \right], \quad {\rm and}\\
& f_a : T_{1} \to \bc \setminus \left( (-\infty,-a] \cup [0,+\infty ) \right), 
\end{split}
\end{equation*}
are one-to-one maps. Define $g_a^0=f_a^{-1}|T_0$ and $g_a^1=f_a^{-1}|T_1$ the corresponding inverse branches of $f_a$ taking values in $T_0$ and $T_1$, respectively.  As for the refinement of $T_0\cup T_1$, we denote by $g_a^{j_1}$ and $g_a^{j_2}$ the appropiate restrictions of $g_a^j$ mapping into $T_{j_1}$ and $T_{j_2}$, respectively.

Assume $z$ is a point of the Julia set whose orbit is entirely contained  in $\cup_{j\in \bbZ} T_j$. We can naturally associate to $z$ the \emph{itinerary} $s(z)=\left(s_{0},s_{1},\ldots \right)$, with $s_{j}\in \bbZ$, if and only if $f_{a}^j(z)\in T_{s_j}$.  Let us concentrate on the space $\Sigma_B=\{0,1\}^\bbN$ of {\it binary sequences} (in what follows \emph{$B$-sequences}). With respect to the refinement of $T_0\cup T_1$, consider the space of {\it extended sequences} given by $\Sigma_E=\{0_{1},0_{2},1_{1},1_{2}\}^\bbN$. Since $f_a$ is a one-to-one map in $T_{0}$ and $T_{1}$, not all extended sequences are {\it allowable}, that is, $f_a$ behaves as a subshift of finite type over the set of points with full orbits inside $T_{0_1}\cup T_{0_2}\cup T_{1_1}\cup T_{1_2}$. Its transition matrix is given by
$$
A=\left(\begin{array}{cccc}1 & 0 & 1 & 0 \\0 & 1 & 0 &1 \\0 & 1 & 0 & 1 \\1 & 0 & 1 & 0 \end{array}\right),
$$
and determines the space of \emph{allowable extended sequences} (in what follows \emph{$A$-sequences}) given by
$$\Sigma_A=\{(s_0,s_1,\ldots)\in \Sigma_E ~|~s_i\in \{0_1,0_2,1_1,1_2\},a_{s_i s_{i+1}}=1, \forall i\}.$$

Denote by $\pi:\Sigma_A\to \Sigma_B$ the \emph{projection map} that transforms an $A$-sequence into a $B$-sequence by erasing all subscripts.  The form of the matrix $A$ makes evident that $\pi$ is a 2-to-1 map. For a given $B$-sequence $t$, denote by $t^1$ and $t^2$ the unique $A$-sequences so that $\pi(t^j)=t$. Observe that by interchanging all subscripts in $t^1$ we obtain $t^2$, and conversely.

\begin{Remark}\label{rem:undefinedAseq}
It is important to observe that points on the real line (or on any of its preimages), do not have well defined $A$-sequences. However, since $\bbR$ is forward invariant under $f_a$, its dynamics and combinatorics are completely understood. Based on the next result, from now on we will only consider $B$-sequences (and its two associated $A$-sequences) that do not end in all zeros. 
\end{Remark}

\begin{Lemma}\label{lem:ends-in-zeros}
Assume $a> 1$ and let $w\in \J(f_a)$ such that $f^k_a(w)\in T_0\cup T_1$ for all $k\geq 0$. Then, $w\in \bbR\cap T_j$ if and only if $s(w)=(j,0,0,\ldots)$, for $j=0,1$.
\end{Lemma}
\begin{proof}
The first implication follows easily by analyzing the action of $f_a$ in $\bbR$. Whenever $a>1$, the set $\A^*(-a)$ intersects the real line in an open interval $(q_a,p_a)$, where $p_a$ is a repelling fixed point and $q_a$ is its only preimage in $\bbR$. Moreover $(-\infty, q_a)\cup (p_a,+\infty)$ consists of points that escape to +$\infty$ along $\bbR^+$ and hence, belong to $\J(f_a)$. Since $f_a$ sends $(-\infty, q_a]$ onto $[p_a,+\infty)$ and this second interval is fixed by $f_a$, then $w$ has a well defined itinerary given by $s(w)=(1,0,0,\ldots)$, if $w\in (-\infty, q_a]\subset T_1$, or $s(w)=(0,0,\ldots)$ if $w\in [p_a,+\infty)\subset T_0$.

To see the second implication, it is enough to show the interval $[p_a,+\infty)$ represents the only set of points  in the Julia set that remain inside $T_0$ for all positive iterates. To do so, we analyze the preimages of $\eta_1$ inside $T_{0_2}$ (the case $\eta_{-1}$ and $T_{0_1}$ is analogous). Since $f_a$ maps $T_{0_2}$ onto the upper half plane Im$(z)>0$, then for each $k\geq 1$, the $k^{\rm th}$ preimage of $\eta_1$ in $T_{0_2}$, namely $\eta_1^k=(g_a^{0_2})^k(\eta_1)$, lies completely inside $T_{0_2}$ (except for its endpoint in $z=-a$) and extends towards infinity into the right half plane. In particular, it lies in the strip bounded by $[-a,+\infty)$ and $\eta_1^{k-1}$ (from bottom to top). Also, note $\eta_1^k$ and $\eta_1^j$ meet only at $z=-a$ whenever $k\neq j$. We claim that $\eta_1^k$ accumulates onto $[p_a,+\infty)$ as $k\to \infty$. For otherwise, we can find a point $x\in [-a,+\infty)$ and $\varepsilon>0$ so $B_\varepsilon(x)\cap \eta_1^k =\emptyset$ for all $k\geq 1$. Nevertheless, since $x$ belongs to the Julia set, it follows from Montel's Theorem the existence of an integer $N>0$ for which $f_a^N(B_\varepsilon(x))\cap \eta_1\neq \emptyset$. Hence, $B_\varepsilon(x) \cap \eta_1^N\neq \emptyset$, a contradiction.

Finally, for any given point $w\in T_{0_2}$ such that $s(w)=(0,0,\dots)$, there exists an integer $m>0$ for which, either $w\in \eta_1^m$ or $w$ lies in the interior of the strip bounded by $\eta_1^{m+1}$ and $\eta_1^{m}$ (from bottom to top). In both situations, $f_a^{m+1}(w)$ lies outside $T_0$. This finishes the proof.
\end{proof}

The rest of this section is devoted to a combinatorial description of the dynamics of points with forward orbits contained in $T_0\cup T_1$ using $A$- and $B$-sequences. First, we focus our study on the set $\A^*(-a)$ and then analyze points in the Julia set that lie far to the right in $T_0\cup T_1$. 
Using known results in complex dynamics we will prove that points with forward orbits completely contained in a given right hand plane are organized into continuous curves and their combinatorics are governed by the transition matrix $A$.
For future reference, we compute the image of a vertical segment bounded above and below by $\zeta_1$ and $\zeta_{-1}$, respectively.

\begin{Lemma} \label{lemma:vertical_segment}
Let $x\in \bbR$ be fixed and consider the vertical segment  $L[x]=\{x+iy \, | \, \zeta_{-1}(x)\leq y \leq \zeta_1(x)\}$. Then $f_a(L[x])$ lies inside the closed round annulus  $f_a(x) \leq |z| \leq f_a(x+\zeta_1(x))$.
\end{Lemma}

\begin{proof}
From the definition of the map $f_a$ we have
\[
|f_a(x+iy)|= a \exp(x+a) \sqrt{(x-(1-a))^2 + y^2}.
\]
\noindent Evidently, when restricted to $L[x]$ for a fixed $x$, the above expression reaches its minimum value when $y=0$ while its maximum value is reached whenever $y=\zeta_1(x)=\zeta_{-1}(x)$. 
\end{proof}

\subsection{Dynamics near $z=-a$}\label{subsection:pol_like}

In \cite{Moro} it was shown that for $a>1$, each Fatou domain of $f_a$ is a bounded, connected component of $\mathcal{A}(-a)$ whose boundary is a Jordan curve. Here we show that $\overline {\mathcal{A}^*(-a)}$ is in fact a quasiconformal image of the closed unit disk.
Precisely, we describe a set of points with bounded orbits inside $T_{0}\cup T_{1}$ through a polynomial-like construction (see \cite{DH}) around the unique and simple critical point $z=-a$.~For technical reasons, we restrict to parameters $a\geq 3$ from now on.

\begin{Proposition} \label{proposition:pol_like}
For any $a\geq3$, there exist open, bounded and simply connected domains $U_a$ and $V_a$ with $ -a \in \overline{U}_a \subset V_a$, such that $(f_{a},U_{a},V_{a})$ is a quadratic-like mapping. Furthermore, the filled Julia set of $(f_{a},U_{a},V_{a})$ is the image under a quasiconformal mapping of the closed unit disk and coincides with $\overline{\mathcal{A}^*(-a)}$.
\end{Proposition}

\begin{figure}[hbt]
  
  \psfrag{a}[][]{\footnotesize $\zeta_1$} 
  \psfrag{b}[][]{\footnotesize $\eta_1$} 
  \psfrag{d}[][]{\footnotesize $\zeta_{-1}$} 
  \psfrag{c}[][]{\footnotesize $\eta_{-1}$} 
  \psfrag{e}[][]{\footnotesize $f_a$} 
  \psfrag{f}[][]{\footnotesize $U_a$}
  \psfrag{g}[][]{\footnotesize $V_a$} 
  \psfrag{h}[][]{\footnotesize $L$}
  \psfrag{i}[][]{\footnotesize $R$}
  \centerline{\includegraphics[width=0.5\textwidth]{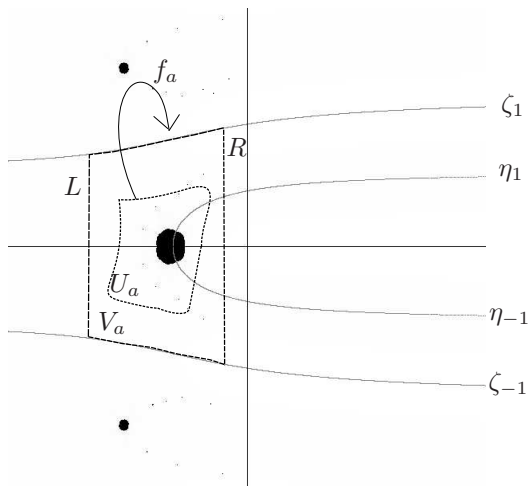}}
   \caption{\small{A sketch of the domains $U_a$ and $V_a$ found in Proposition \ref{proposition:pol_like} }}
      \label{fig:polynomial}
\end{figure}

\begin{proof}
Define $V_a$ as the open, simply connected pseudo-rectangle given by
$$
V_{a}=\{z \in \bc \ | \    -a-6\ln a < {\rm Re}(z) < \frac{1-a}{2}, \ \zeta_{-1}({\rm Re}(z)) < {\rm Im}(z) < \zeta_{1}({\rm Re}(z)) \}.
$$
First, we show that $V_{a}$ maps outside itself. Indeed, the top and bottom boundaries of $V_a$ map into a segment lying on the positive real line, thus outside $V_a$ as $(1-a)/2\leq -1$. Also, note that $V_a$ lies in the interior of the annulus
$$\frac{|1-a|}{2}\leq |z| \leq |-a-6 \ln a+ i2\pi|.$$
Following the notation in Lemma \ref{lemma:vertical_segment}, $L=L[-a-6\ln a]$ and $R=L[(1-a)/2]$ are the left and right hand boundaries of $V_a$. We show next the images of $L$ and $R$ lie in the complementary components of the annulus. First, for $z\in L$ we have
$$
|f_{a}(z)|=\frac{1}{a^5}\sqrt{\left( 1+6\ln a  \right)^2+y^2} < \frac{1}{a^5}\sqrt{\left( 1+6\ln a  \right)^2+4\pi^2}<\frac{1}{2} < \frac{|1-a|}{2}.
$$

Similarly, if $z\in R$
\[
|f_{a}(z)|=a e^{\frac{a+1}{2}}\sqrt{\left(\frac{1-a}{2}  \right)^2 + y^2}>a e^{\frac{a+1}{2}} \frac{|1-a|}{2}  >\sqrt{(a+6 \ln a)^2 + 4\pi^2} 
\]
\noindent for all $a\geq 3$, proving thus that $V_a$ is mapped outside itself under $f_a$.

Now, we define $U_a$ to be the connected component of $f^{-1}_a(V_a)$ containing $-a$. Since $-a$ is a superattracting fixed point with multiplicity one, and  there are no other critical points, it follows that $\overline{U_a}\subset V_a$ and the map $f_a:U_a \to V_a$ sends $\partial U_a$ to $\partial  V_a$ with degree 2, as $-a$ is a simple critical point. We conclude that $(f_a,U_a,V_a)$ is a quadratic-like mapping.

What is left to verify is that the filled Julia set of $(f_a,U_a,V_a)$ is a quasi-disk. Recall that the filled Julia set of a polynomial-like mapping is defined as the set
$\{z\in U_a~|~f_a^n(z)\in U_a~\text{for~all}~n\geq 0\}$.
Being $(f_a,U_a,V_a)$ a quadratic-like mapping, there exists a quasiconformal conjugacy with a polynomial of degree two that has a superattracting fixed point.~Thus, the polynomial must be $z\mapsto z^2$ after a holomorphic change of variables, if necessary.~So the filled Julia set of $(f_a,U_a,V_a)$ is the image under a quasiconformal mapping of the closed unit disk.
\end{proof}

\begin{Proposition}\label{prop:consequences_pol_like}
Let $a\geq3$.   The following statements hold.
\begin{enumerate}
\item[(a)] The map $f_{a}$ restricted to the boundary of $\mathcal A^{*}(-a)$ is conjugate to the map $\theta \mapsto 2\theta$ in the unit circle.
\item[(b)] Let $t\in \Sigma_E$ be an extended sequence that does not end in $0_i$'s. Then, $t$ is an $A$-sequence if and only if  there exists a unique point $z \in \partial \mathcal A^{*}(-a)$ that realizes $t$ as its itinerary.
\end{enumerate}
     
\end{Proposition}

\begin{proof}
Statement (a) is a direct consequence of the previous proposition since the map $f_a$ is conjugate in $\partial \mathcal A^*(-a)$ to $z \mapsto z^2$ acting on the unit circle. 

We prove statement (b) by defining a partition of the boundary of $\mathcal A^*(-a)$ that coincides with the refinement of the partition $T_0 \cup T_1$ discussed before. For simplicity, angles are measured by $[0,1]$. Denote by $z(0)$, $z(1/4)$, $z(1/2)$ and $z(3/4)$ the points in $\partial \mathcal A^*(-a)$ corresponding under the conjugacy between $f_a$ and $z \mapsto z^2$ to points in $S^1$ of angle $\theta=0,1/4,1/2 \,$ and $\,3/4$. Now label points in $\partial \mathcal A^*(-a)$ in the following way: traveling along $\partial \mathcal A^*(-a)$ in a clockwise direction, associate the symbol $0_1$ to the arc joining $z(0)$ and $z(3/4)$, the symbol $1_1$ to the arc joining $z(3/4)$ and $z(1/2)$, $1_2$ to the arc joining $z(1/2)$  and $z(1/4)$, and $0_2$ to the arc joining $z(1/4)$ and $z(0)$. We leave it to the reader to verify the transition matrix for this partition under the action of $f_a$ is exactly $A$ and the labeling is consistent with the one defined by the $T_{j_i}$'s.
\end{proof}

\begin{Theorem}\label{thm:BdOrbit}
Let $z$ be a point such that $f_a^n(z)\in \overline{T_0\cup T_1}$ for all $n\geq 0$. If $z$ belongs to the Julia set and has bounded orbit, then $z\in \partial \A^*(-a)$.
\end{Theorem}

\begin{proof}
If $z\in \J(f_a)$ satisfies the hypotheses, we can find $m<0<M$ so that $m< \text{Re}(f_a^n(z)) < M$ for all $n\geq 0$.

Let $\varepsilon>0$ small enough and denote by $B_\varepsilon=\overline{B_\varepsilon(0)}$. Since the orbit of the origin escapes monotonically along the positive real line, redefining $M$ if necessary, there exists an integer $N=N(M)>0$ for which $B_\varepsilon, f_a(B_\varepsilon),\ldots, f^{N}_a(B_\varepsilon)$ are pairwise disjoint compact domains, such that for all $0\leq j \leq N-1$,
\[
\begin{array}{ll}
& f^j_a(B_{\varepsilon}) \subset T_0 \cap \{ z \, |\, {\rm Re}(z)<M\},  \hbox{ and }\\
& f^N_a(B_{\varepsilon}) \subset T_0 \cap \{ z \, | \, {\rm Re}(z)> M\}.
\end{array}
\]
Moreover, we may choose $\varepsilon$ small enough so $B_\varepsilon, f_a(B_\varepsilon),\ldots, f^{N}_a(B_\varepsilon)$ are all compact domains contained in $T_0$.
Finally, select $m'\leq m<0$ so the subset $\{z\in T_1~|~\text{Re}(z)\leq m'\}$ maps completely inside $B_\varepsilon$.

Denote by $\Phi_a:\bbD\to \A^*(-a)$ the B\"ottcher coordinates tangent to the identity at the origin. For $0<r<1$ let $\Delta_r=\Phi_a(B_{r}(0))$. 
Clearly $\Delta_r\subset \A^*(-a)$ and maps compactly into its own interior. Moreover, we can choose $r$ small enough so for $j=0,1$, $T_j \setminus \Delta_r$ is a connected set and the intersection of $\Delta_r$ and $\bbR$ is an open interval $(c,d)$, since $\Phi_a$ has been chosen to be tangent to the identity at the origin. We can now define the set $E$, illustrated in Figure~\ref{fig:periodic}, as follows
\[ E=\{z\in \bbC~|~m'< \text{Re}(z)< M,~\zeta_{-1}(\text{Re}(z))< \text{Im}(z)< \zeta_1(\text{Re}(z))\} \setminus \Omega  
\]
\noindent where $\Omega =\bigcup_{k=0}^{N-1} f^k_a(B_{\varepsilon})\cup \Delta_r \cup [c,M]$.  It is easy to verify that $E$ is an open, bounded, connected and simply connected set and $\partial  \A^*(-a) \subset E$. Moreover, $E\subset f_a(E)$ although some boundary components map into $\partial E$. Indeed, $f_a(\partial E)\cap \partial E$ consists of segments along the real line and $f^j_a(\partial B_{\varepsilon})$, for $j=1,2,\ldots,N-1$. So after $N$ iterations, the only boundary points mapping into $\partial E$ are points over the real line, thus having $B$-itineraries $(1,0,\ldots)$ or $(0, 0,\ldots)$. We show next that for $\ell>0$ sufficiently large, $f_a^{-\ell}|E$ becomes an strict contraction.

\begin{figure}[hbt]
  
  \psfrag{a}[][]{\footnotesize $\partial A^*(-a)$} 
  \psfrag{r}[][]{\footnotesize $-a$} 
  \psfrag{n}[][]{\footnotesize $\Delta_r$} 
  \psfrag{h}[][]{\footnotesize $M$} 
  \psfrag{b}[][]{\footnotesize $-a$} 
  \psfrag{l}[][]{\footnotesize $B_{\varepsilon}(0)$} 
  \psfrag{c}[][]{\footnotesize $E$} 
  \psfrag{d}[][]{\footnotesize $m'$} 
  \centerline{\includegraphics[width=0.5\textwidth]{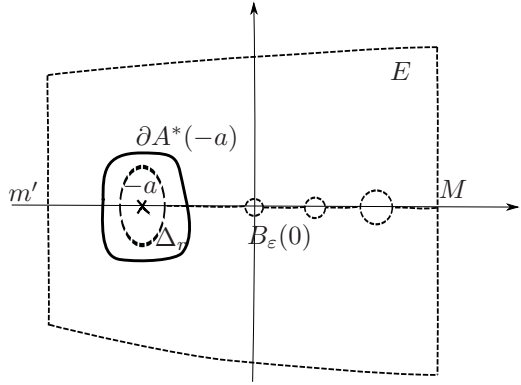}}
   \caption{\small{An schematic representation of the set $E$ described in Theorem~\ref{thm:BdOrbit}.}}
      \label{fig:periodic}
\end{figure}

Let $s=(s_0, s_1, s_2,\ldots)$ be the $A$-sequence associated to $z$. If $s$ has finitely many $1_i$'s, then its $B$-sequence ends with $0$'s and since $z$ has bounded orbit inside $\overline{T_0\cup T_1}$, then $z\in \partial \A^*(-a)$ by Lemma~\ref{lem:ends-in-zeros}. 
If $s$ has infinitely many $1_i$'s there exists a first integer $n>N$ for which $s_{n}\in\{1_1,1_2\}$, and thus $\overline{E\cap T_{s_{n}}} \subset f^{n}(E)$, as points in $f^{n}(\partial E)$ mapping into $\partial E$ have by now itinerary $(0,0,\ldots)$.

For each $k\in \{0_1, 0_2, 1_1, 1_2\}$, the set $E\cap T_k$ is an open, connected and simply connected set with a Riemann mapping given by $\psi_k:\bbD\to E\cap T_k$. Consider the mapping $\Psi_{\ell}:\bbD\to \bbD$ with $\ell>n$, given by
$$\Psi_{\ell} = \psi_{s_0}^{-1}\circ( g_a^{s_0}\circ \ldots \circ g_a^{s_{\ell-1}} )\circ \psi_{s_\ell}.$$

It follows that $\Psi_\ell(\bbD)$ is compactly contained in $\bbD$, that is, for $\ell>n$, $\Psi_\ell$ is an strict contraction with respect to the Poincar\'e metric on the unit disk. Consequently, the sets $\overline{\Psi_\ell(\bbD)}$ form a nested sequence of compact sets with diameters converging to zero as $\ell\to \infty$. This implies that for any $w\in \bbD$, $\lim_{\ell\to \infty} \Psi_\ell(w)$ exists and is independent of the point $w$. 
Therefore, by construction, $\psi_{s_0}(\lim_{\ell\to \infty} \Psi_\ell(0))=z$ is the unique point in $E$ with itinerary $s$. From Proposition~\ref{prop:consequences_pol_like}(b) we conclude $z\in \partial \A^*(-a)$.
\end{proof}

\subsection{Dynamics near infinity}
\label{subsection:tails_to_the_right}

Our first aim is to prove that for $R>0$ sufficiently large and the region
$$H_R = \{z\in T_0\cup T_1~|~ {\rm Re}(z)\geq R\},$$
there exist continuous curves in $\J(f_a)\cap H_R$ consisting of points whose orbits escape to $+\infty$ with increasing real part. These curves are usually known as \emph{tails}. The existence of tails as disjoint components of the Julia set were first observed by Devaney and Tangerman \cite{DT}  for certain entire transcendental maps and by Schleicher and Zimmer \cite{SZ} for the exponential family $E_\lambda(z)=\lambda \exp(z)$ and all $\lambda\in \bbC$. In greater generality, Bara\'nski \cite{Ba} and Rempe \cite{R} have shown the existence of tails for hyperbolic maps belonging to the class $\cal B$. 

For completeness, we analyze in detail some of their results in the setting of our work to obtain tails in $\J(f_a)\cap H_R$. Once each tail has been assigned an $A$-sequence, we describe a pullback process to compute the full set of points in $T_0\cup T_1$ associated to such $A$-sequence. In the final section, we study the topological properties of that set.

Consider an entire transcendental map $f$ in the class $\cal B$. 
The \emph{escaping set of $f$}, denoted as $I(f)$, is the set of points whose orbits under $f$ tend to infinity. For an entire transcendental map, Emerenko \cite{E} has shown the Julia set coincides with the boundary of $I(f)$.
We say that two maps $f,g\in \cal B$ are \emph{quasiconformally equivalent near infinity} if there exist quasiconformal maps $\phi_1, \phi_2:\bbC \to \bbC$ that satisfy $\phi_1 \circ f=g\circ \phi_2$ in a neighborhood of infinity.

\begin{Theorem}[Rempe, 2009]\label{theorem:lasse}
Let $f,g \in \mathcal B$ be two entire transcendental maps which are quasiconformally equivalent near infinity. Then there exist $\rho>0$ and a quasiconformal map $\theta:\bc \to \bc$ such that $\theta \circ f = g \circ \theta$ on
$$
A_{\rho}=\{ z\in\bc \ | \ |f^n(z)|>\rho, \ \forall n \geq 1   \}.
$$
Furthermore, the complex dilatation  of $\theta$ on $I(f)\cap A_{\rho}$ is zero.
\end{Theorem}

A straightforward computation shows that $f_{a}(z)$ is conjugate to the function $\tilde{f}_{a}(z)=az\exp(z+1)-(1-a)$ under the conformal isomorphism $\varphi_a(z)=z-(1-a)$. Since $\tilde{f_a}(z)= \varphi_a \circ f_a \circ \varphi_a^{-1}(z)$, it is easily verified that $\tilde{f}_{a}$ has a free asymptotical value at $z=a-1$ and a fixed critical point at $z=-1$. 
In turn, $\tilde{f}_{a}$ is (globally) conformally equivalent to $g_{b}(z)=bz\exp(z)$ via $\phi_{1}(z)=\alpha z+\alpha(1-a)$ and $\phi_{2}(z)=z$. Indeed, it is easy to see that $\phi_{1}\circ \tilde{f}_{a} = g_{b}\circ \phi_{2}$, where $b = e a\alpha$.

For small values of $b$, the Fatou set of $g_{b}$ consists solely of the completely invariant basin of attraction of the fixed point (and asymptotic value) $z=0$. Thus, we can describe the Julia set of $g_b$ by applying the following result found in \cite{Ba}.

\begin{Theorem}[Bara\'nski, 2007]\label{theo:baranski}
Let $g$ be an entire transcendental function of finite order so that all critical and asymptotic values are contained in a compact subset of a completely invariant attracting basin of a fixed point. Then ${\cal J}(g)$ consists of disjoint curves (hairs) homeomorphic to the half-line $[0,+\infty)$. Moreover, the hairs without endpoints are contained in the escaping set $I(g)$.
\end{Theorem}

These disjoint curves are usually known as hairs or \emph{dynamic rays}. Several consequences are derived from the above theorem. Firstly, if $\gamma$ denotes a hair, it can 
be parametrized by a continuous function $h(t), \ t\in [0,+\infty)$, such that $\gamma=h([0,+\infty))$. The point $h(0)$ is called the \emph{endpoint} of the hair. Secondly, each hair is a curve that extends to infinity and, for $r>0$, we say that  $\omega=h((r,+\infty))$ is the {\it tail of the hair}. Moreover, all points in a given hair share the same symbolic itinerary defined by a dynamical partition of the plane with respect to $f$, and for every point $z\in \gamma$ that is not the endpoint, we have $f^n(z)\to \infty$ as $n\to \infty$. Finally, if $f$ and $g$ are as in Theorem  \ref{theorem:lasse} and in addition, $g$ satisfies hypotheses in Theorem \ref{theo:baranski}, then near infinity the topological structure of the escaping set of $f$ is also given by disjoint curves extending to infinity. We refer to \cite{BJR} for
a topological description of the Julia set in terms of what is known as Cantor bouquets. Furthermore, the dynamics of $f$ in those curves is quasiconformally conjugate to the dynamics of $g$ in the corresponding curves near infinity. We deduce the following result based on the previous theorems and the specific expression of $f_a$.

\begin{Proposition}\label{prop:ConjugTails}
Let $f_a(z) = a(z-(1-a))\exp(z+a)$, $g_b(z)=bz\exp(z)$, $a\geq 3$, and $b$ a complex parameter.
\begin{enumerate}
\item[(a)] If $|b|$ is small enough, the Julia set of $g_{b}$ is given by the union of  disjoint hairs. Each hair lands at a distinguished endpoint. Moreover, hairs without endpoints are contained in $I(g_b)$. 
\item[(b)] Let $R>0$ large enough. The set of points with forward $f_a$-orbits that are always contained in $H_R$ are given by the union of disjoint curves extending to infinity to the right. All points in those curves belong to $I(f_a)$.
Precisely, these curves are quasiconformal copies of connected components of hairs described in (a).
\item[(c)] To each curve in (b) that does not coincide with  $\bbR$ or one of its  preimages, we can assign a unique sequence  $t$ in $\Sigma_A$. We denote this curve by $\omega_t$. All points in $\omega_t$ escape to infinity under the action of $f_a$ following the  itinerary $t$.
\item[(d)]  For any $t\in \Sigma_A$, if $R>0$ is large enough, there exists a unique curve $\omega_t$ in $H_R$. Moreover, for each $r \geq R$, $\omega_t \cap \{z~|~\text{Re}(z)=r\}$ is a unique point. In particular $\omega_t$ is the graph of a function.
\item[(e)] $\omega_t$ is a tail, i.e., a quasiconformal copy of a tail in (a).
\end{enumerate} 
\end{Proposition}

\begin{figure}[hbt]  
  \psfrag{a}[][]{\footnotesize $\zeta_1$} 
  \psfrag{b}[][]{\footnotesize $\zeta_{-1}$} 
  \psfrag{g}[][]{\footnotesize $\eta_{-1}$} 
  \psfrag{l}[][]{\footnotesize $\eta_{1}$} 
  \psfrag{h}[][]{\footnotesize {\rm Re}$(z)=R$}  
  \psfrag{w}[][]{\footnotesize $w_t$}
  \psfrag{f}[][]{\footnotesize $\alpha_t$}
  \psfrag{c}[][]{\footnotesize $|z|=f_a(R+i\zeta_1(R))$}
  
  \centerline{\includegraphics[width=0.4\textwidth]{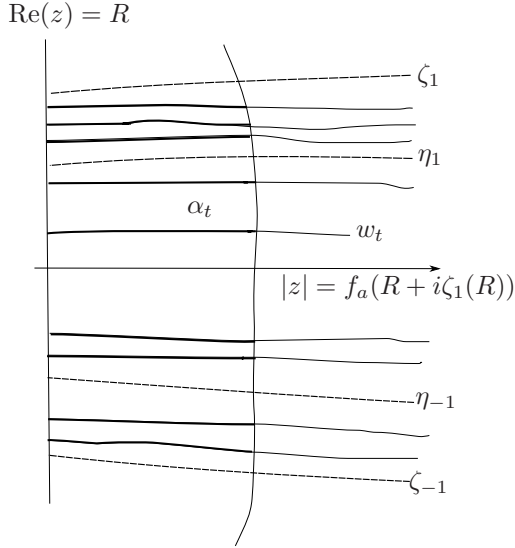}}
   \caption{\small{The set $H_R$, the tail $\omega_t$ and its base $\alpha_t$ for some $t\in \Sigma_A$. Bases are depicted as bold lines. }}
      \label{fig:tails}
\end{figure}

\begin{proof}

Statement (a) follows directly from Theorem~\ref{theo:baranski}. To see statement (b) note that for any value of $a\neq 0$ and $b$, there exists $\alpha=b/(ea)$ so $f_a$ and $g_b$ are (globally) conformally equivalent. By Theorem \ref{theorem:lasse}, there exists $\rho>0$ such that $f_a$ and $g_b$ are conjugate on the set $A_\rho=\{z~|~|f_a^n(z)|>\rho, \forall n \geq 1\}$. Let $R>\rho$. Denote by $S$ the set of points in $H_R$ with forward $f_a$-orbits contained in $H_R$. Clearly, each point in $S$ must belong to $\J(f_a)$ (since a point in the Fatou set eventually maps into $\A^*(-a)$) and in particular, $S\subset \J(f_a)\cap A_R$. Moreover, far enough to the right, a point $z\in H_R$ whose forward orbit remains forever in $H_R$ must satisfy Re$\left(f_a^{k+1}(z)\right) > {\rm Re}\left(f_a^k(z)\right)$ for all $k>0$ (see Lemma  \ref{lemma:vertical_segment}). Hence from Theorems  \ref{theorem:lasse} and \ref{theo:baranski} we know that  $S$ is the union of disjoint curves extending to infinity (that is, quasiconformal copies of components of the hairs in (a)) belonging to the escaping set.


To prove statement (c) we start by assigning to each of these curves, $\omega$, a unique sequence in $\Sigma_A$. Let $z_0$ be a point in $\omega$ and let $s(z_0)$ be its itinerary  in $\Sigma_A$ in terms of the partition $T_{0_1}, T_{0_2}, T_{1_1}$ and $T_{1_2}$.  By assumption  this itinerary is well defined since $\omega$ is not  $\bbR$ or one of its  preimages. Let $z_1$ be another point in $\omega$.  We claim that $s(z_1)=s(z_0)$ and proceed by contradiction. Let $\mathcal C$ be the connected component of $\omega$ joining $z_0$ and $z_1$. If $s(z_1)\neq s(z_0)$, there exists an integer $k\geq 0$ for which the $k^{\text th}$ entries in both itineraries are the first ones to differ.

Hence there is a point $q\in \mathcal C$ such that $f_a^k(q)$ belongs to  either $[R,\infty)$, $\eta_1$ or $\eta_{-1}$.  Clearly $f_a^k(q)$ cannot belong to $\eta_{\pm 1}$, since otherwise $f_a^{k+1}(q)\in \bbR^-$ and by hypothesis the forward orbit of $q$ belongs  to $H_R$. On the other hand $f_a^k(q)$ cannot belong to $\bbR$ since by item (b), $f_a^k\left(\omega\right)$ and $\bbR$ are disjoint curves of the escaping set. 

Finally, there cannot be two curves having the same itinerary as this will imply the existence of an open set of points following the same itinerary, which is impossible. Thus, we may now denote by $\omega_t$ the unique curve in $H_R$ formed by escaping points with itinerary $t$.

To  prove statement (d) fix  $t=(t_0,t_1,\ldots, t_n,\ldots) \in \Sigma_A$ and $R$  large enough. Let $r\geq R$ and let $L[r]=\{r+iy \, | \, \zeta_{-1}(r)\leq y \leq \zeta_1(r)\}$. We denote by $I_{t_0}$ the set of points in $L[r] \cap \overline{T_{s_0}}$.  We know that the image of the vertical segment $L[r]$ cuts across $H_R$ in two almost vertical lines (see Lemma \ref{lemma:vertical_segment}). Recall that when $f_a$ is restricted over the set of points with forward orbits in $H_r$, it behaves as a subshift of finite type governed by the matrix $A$, and when restricted to $T_0$ or $T_1$ it is a one-to-one map. Using these facts, we can find a unique subinterval $I_{t_0t_1}\subset I_{t_0}$ formed by points with forward orbits inside $H_r$ and itineraries starting as $(t_0,t_1,\ldots)$. Inductively, for each $n>0$, $I_{t_0t_1\ldots t_n}$ is the unique subinterval of $I_{t_0t_1\ldots t_{n-1}}$ formed by points with forward orbits inside $H_r$ and itineraries starting with $(t_0,t_1,\ldots, t_n,\ldots )$. Clearly,
$$
I_{t_0t_1,\ldots t_n} \subset I_{t_0t_1\ldots t_{n-1}} \subset \ldots \subset  I_{t_0t_1} \subset I_{t_0}.
$$
Due to the expansivity of $f_a$ in $H_R$ we obtain $\cap_{n\geq 0} I_{t_0t_1,\ldots t_n} =\{q\}$, and by construction $q$ must have itinerary $t$. Moreover $f^k(q)\to \infty$ as $k\to \infty$, so $q$ must belong to the unique curve $\omega_t$ described in (c).  The above arguments imply that $\omega_t$ intersects Re$(z)=r$ at a unique point. Therefore we can parametrize $\omega_t$ on the interval $[r,\infty)$.

Finally to see statement  (e) we observe that there are no endpoints associated to bounded itineraries in $H_R$, implying that each $w_t$ is a quasiconformal copy of the tail of some hair in (a). On one hand, the only points in $T_0 \cup T_1$ with bounded orbit belong to $\partial A^*(-a)$, so they are not in $H_R$. On the other hand if there were and endpoint {\it of a hair} with orbit escaping to infinity in $H_R$, it should have an itinerary, say, $t$. But from statement (d)  there is a (unique) curve $w_t$ going from Re$(z)=R$ to infinity with such itinerary $t$, a contradiction. 

\end{proof}
\begin{Definition}
Given any $A$-sequence $t$, each curve described in Proposition~\ref{prop:ConjugTails} will be denote by $\omega_{t}=\omega_{t}(R)$ and it will be called the \emph{tail with itinerary $t$} contained in the half plane Re$(z)\geq R$.  

The component of $\omega_t$ that cuts across
$$F_R=\{z\in H_R~|~|z|<f_a(R+i\zeta_1(R))\}$$
is called the \emph{base} of the tail, and it will be denote by $\alpha_t=\alpha_t(R)$.
\end{Definition}
See Figure \ref{fig:tails}. 


We now describe a pullback construction to extend the tail $\omega_t$ into a longer curve. Recall $g_a^{j_i}=f_a^{-1}|T_{j_i}, j_i\in\{0_1,0_2,1_1,1_2\}$ and consider the shift map $\sigma: \Sigma_A\to \Sigma_A$ acting on the space of $A$-sequences. Let $t=(t_0,t_1,\ldots)\in \Sigma_A$. By the conjugacy of $f_a | \bigcup_{t\in \Sigma_A} \omega_t$ with $\sigma|\Sigma_A$, it follows that $\omega_{\sigma (t)}$ properly contains $f_{a}(\omega_t)$ since this curve lies in $H_R\setminus F_R$ (see Proposition \ref{prop:ConjugTails}). Hence $f_a(\omega_t)$ is a curve that misses the base $\alpha_{\sigma(t)}$. Consequently, $g_{a}^{t_0}(\omega_{\sigma(t)})$ is a continuous curve that lies in $T_{t_0}$ and extends $\omega_t$ to the left of Re$(z)=R$. Clearly, any point in the extended curve $g_{a}^{t_0}(\omega_{\sigma (t)})$ has itinerary $t$. Inductively, consider
 $$
 g_{a}^{t_0}\circ \ldots \circ g_{a}^{t_{n-1}}\left(\omega_{\sigma^n(t)}\right).
 $$
This pullback iteration is always defined as long as the extended curve does not meet $z=0$, which is impossible since $\br$ is forward invariant. We thus obtain a curve of points with itinerary $t$, and each pullback iteration extends its predecessor. 

\begin{Definition}
Let
\begin{equation}\label{eq:pullback-hair}
\gamma (t)=\bigcup_{n=0}^\infty g_{a}^{t_0}\circ \ldots \circ g_{a}^{t_{n-1}}\left(\omega_{\sigma^n(t)}\right).
\end{equation}
We call $\gamma (t)$ the {\it hair\/} associated to $t$.
\end{Definition}
Note that the pullback process described before may or may not produce an endpoint. In Section~4 we show that in some cases, $\gamma(t)$ is a hair with an endpoint in $\partial \A^*(-a)$, and in some other cases $\gamma(t)$ is a non-landing hair and accumulates everywhere upon itself. By the following theorem found in \cite{C}, we will conclude that $\overline{\gamma(t)}$ is an indecomposable continuum.

\begin{Theorem}[Curry, 1991]\label{theo:curry}
Suppose that $X$ is a one-dimensional nonseparating plane continuum which is the closure of a ray that limits on itself. Then $X$ is either an indecomposable continuum or the union of two indecomposable continua.
\end{Theorem}

\begin{Remark}
A \emph{ray}  is defined as the image of $[0,+\infty)$ under a continuous, one to one map. Giving any positive number $\alpha$, the image of $[\alpha,+\infty)$ under the same map is known as a \emph{final segment} of the ray. Then, the ray \emph{limits on itself} if it is contained in the closure of any final segment of itself.
\end{Remark}

\section{Targets in $H_{R}$}\label{section:targets}

The main result in this section will be Theorem~\ref{theorem:indecom}, where we construct $B$-sequences so their associated $A$-sequences produce hairs that accumulate everywhere on themselves.\\

We now set up \emph{targets} around the $n^{\rm th}$ \emph{image of the base} $\alpha_{t}$ of each tail $\omega_t$. The construction is very similar to the one presented in Devaney and Jarque \cite{DJ2} and Devaney, Jarque and Moreno Rocha in \cite{DJM}, although the existence of a critical point in our present case requires some modifications.\\

Let $t$ be a given $A$-sequence. We first enlarge the base $\alpha_t$ inductively along $\omega_t$. Set $\alpha_{t,0}=\alpha_t$ and consider the two bases $\alpha_{\sigma^{-1}(t)}$.  The set
$f_{a}(\alpha_{\sigma^{-1}(t)})$ (taking the two possible bases) is a subset of $\omega_t$. Thus, the set  
$$
\alpha_{t,1} = \alpha_{t,0} \cup f_{a}(\alpha_{\sigma^{-1}(t)}),
$$
is an extension of $\alpha_{t,0}$ along $\omega_t$.
Inductively, define the $n^{\rm th}$ image of the base $\alpha_t$ as
$$
\alpha_{t,n} = \alpha_{t,n-1} \cup f_{a}^n(\alpha_{\sigma^{-n}(t)}).
$$
It is easy to verify that $\{\alpha_{t,n}\}_{n\in \bn}$ is a sequence of curves satisfying the following three conditions:
\begin{enumerate}
\item[(i)] $\alpha_{t,0} = \alpha_t$,
\item[(ii)] $\alpha_{t,n} \subset \alpha_{t,n+1}$, and
\item[(iii)] $\bigcup\limits_{n\geq 0}\alpha_{t,n} = \omega_t$.
\end{enumerate}

In order to define a target around each $\alpha_{t,n}$, consider $\xi,\eta\in \bbR^+$ and let
\[
V(\xi,\eta) = \{ z \in H_R~|~\xi - 1 <{\rm Re}(z) < \eta+1 \}.
\]
By definition $V( \xi,  \eta)$ is a rectangular region bounded above and below
by components of  $\zeta_{-1}$ and $\zeta_{1}$, respectively.

 \begin{Lemma}\label{lemma:macrolemma}
Let $R>0$ be large enough. For all $n \geq 0$ there exist positive real numbers $\xi_n$ and $\eta_n$ such that the following statements hold.
\begin{enumerate}
\item[(a)]  For every $t\in \Sigma_A$, the $n^{\rm th}$ iterate of $\alpha_t$ belongs to the interior of $V(\xi_n,\eta_n)$.
\item[(b)]  For every $\ell \geq 0$, $V(\xi_{n+1},\eta_{n+\ell+1})$ is compactly contained inside $f_a(V(\xi_n,\eta_{n+\ell}))$.
\end{enumerate} 
\end{Lemma}

\begin{proof}  
Set $\eta_0 = f_a\left(R+i\zeta_{1}(R)\right)$, $\eta_{n+1}=f_a\left(\eta_n+i  \zeta_{1}\left(\eta_{n}\right)\right)$ and $\xi_n=f_a^n(R)$ for $n\geq 0$. It remains to verify that these values satisfy statements (a) and (b). Observe that for $R$ large enough, the image of every vertical segment  $L[R]$ cuts across  $T_{0}\cup T_{1}$ in two almost vertical lines: one {\it near} Re$(z)=f_{a}(R)$ and the other {\it near} Re$(z)=f_{a}\left(R+i\zeta_{1}(R)\right)$ (see Lemma \ref{lemma:vertical_segment}).  

Statement (a) follows directly from the definition of $\xi_n$ and $\eta_n$, since at each step we choose $\xi_n$ and $\eta_{n}$ to be respectively the smallest and largest possible values of Re$(f^n_a(\alpha_t))$ for all $t\in \Sigma_A$, and moreover it is easy to check that $\xi_{n}<\eta_{n-1}$.

We prove statement (b) when $\ell =0$. The case $\ell >0$ follows similarly. The proof of the statement proceeds in two steps. The first one is to verify the inequality 
\[
f_a\left(\xi_n-1+i\zeta_{1}\left(\xi_n-1\right) \right) < \xi_{n+1}-1 \,.
\]
From the definition of $f_a$ we obtain that 
\begin{eqnarray*}
f_a\left(\xi_n-1+i\zeta_{1}\left(\xi_n-1\right) \right) & < &a e^{a+\xi_n-1} \sqrt{(\xi_n-1 -(1-a))^2+4\pi^2}\\
& < & a e^{a+\xi_n} (\xi_n-(1-a))-1\\
& = & \xi_{n+1}-1,
\end{eqnarray*}
where the second inequality is satified if $R$ is large enough. The second step is to verify the inequality 
\[
f_a(\eta_n+1) >  \eta_{n+1}+1 = f_a(\eta_n+i\zeta_1(\eta_n))+1.
\]
By evaluating both sides we obtain
\[
a e^{a+\eta_n+1} (\eta_n+a) > a e^{a+\eta_n}\sqrt{(\eta_n+a-1)^2+\left(\zeta_1(\eta_{n})\right)^2} +1.
\]
Consequently the image of $V(\xi_n,\eta_n)$ contains $\overline{V(\xi_{n+1},\eta_{n+1})}$ as desired.
\end{proof}

\begin{Definition}
The set $V(\xi_n,\eta_n)$ will be called the $n^{\rm th}$ target of $f_a$.
\end{Definition}

Targets provide a useful tool in the proof of Theorem~\ref{theorem:indecom}. Before getting into the details, we briefly explain why targets are so important in our construction. For large values of $n$, an $n^{\rm th}$ target corresponds to a rectangular region with arbitrarily large real part and imaginary part bounded, in absolute value, by $2\pi$. Then, we may pullback the $n^{\text th}$ target using suitable branches of $f_a^{-1}$ to obtain two sequences of nested subsets inside $V(\xi_0,\eta_0)$. We intend to prove each nested sequence contains not only a base $\alpha_{t^{i}}, i=1,2,$ but there are also other components of $\gamma(t^i)$ accumulating into the base.\\

As observed before, each target $V(\xi_n,\eta_{n+\ell})$ intersects the
fundamental domains $T_{0}$ and $T_{1}$ for any $n\in \bz^+$.~Denote by $W_{n,\ell}^{0_{i}}$ and $W_{n,\ell}^{1_{i}}$ the domains given by  $V(\xi_n,\eta_{n+ \ell})\cap T_{0_{i}}$ and $V(\xi_n,\eta_{n+\ell})\cap T_{1_{i}}$, respectively. 
The next step is to show that, by considering appropriate preimages of these $W$-sets, we obtain a nested sequence of neighborhoods around two particular bases of tails.

In what follows, we will work solely with $B$-sequences and its respective $A$-sequences. We may also assume that every $B$-sequence has infinitely many 1's to avoid taking a preimage of the positive real line.
The next result is a suitable restatement of Lemma 4.3 in \cite{DJM} following our notation. 

\begin{Lemma}\label{lemma:two-nested}
Let $\ell \ge 0$ and let $t=(\tau_0,1,\tau_1,1,\tau_2,1,\ldots)$ be a $B$-sequence, where each $\tau_k=\left(\tau_k^1,\ldots, \tau_k^{n_k}\right)$ denotes a finite block of binary symbols $\{0,1\}$ of length $n_k$. Let $m_j=j+1+n_0+n_1+\ldots + n_j, \ j\geq 0$.  Then for each  $m_{j}$ the sets
\begin{equation}\label{eq:W-preimage-usingB}
g_{a}^{\tau_{0}^1}\circ\cdots\circ
g_{a}^{\tau_{j}^{n_j}} \ \left(W_{m_{j},\ell}^{1_i} \right), \ i=1,2,
\end{equation}
form two nested sequences of subsets of $V(\xi_0,\eta_\ell)$. Moreover, if we denote by  $t^{1}$ and $t^{2}$ the two allowable $A$-sequences satisfying $\pi\left(t^{1}\right)=\pi \left(t^{2}\right)=t$, then $\alpha_{t^1,\ell}$ and $\alpha_{t^2,\ell}$ are contained each in a nested sequence of subsets of $V(\xi_0,\eta_\ell)$ given by (\ref{eq:W-preimage-usingB}).
\end{Lemma}

\begin{proof} 
Due to Lemma \ref{lemma:macrolemma}, for each $m_j$ and  $i=1,2$, the sets  $W_{m_{j},\ell}^{1_i} \subset V(\xi_{m_j},\eta_{m_j+\ell})$, can be pulled back following $(\tau_0,1,\tau_1,1,\ldots, \tau_j)$. These preimages yield two nested sequences of subsets in  $V(\xi_0,\eta_\ell)$, one sequence corresponds to  preimages of $W_{m_{j},\ell}^{1_1}$ and the other to  preimages of $W_{m_{j},\ell}^{1_2}$. Notice there is a unique way of pulling back each $W_{m_{j},\ell}^{1_i},\ i=1,2$ following the $B$-sequence $t$, since ${f_{a}}|{T_k},\ k=0,1$ are one-to-one maps.  
Also, at each step  of the construction, the points belonging to these nested subsets are points of $V(\xi_0,\eta_\ell)$  with $B$-itinerary $t=(\tau_0,1,\tau_1,1,\ldots, \tau_j,\ldots)$. So, $\alpha_{t^1,\ell}$ and $\alpha_{t^2,\ell}$ must be inside all of them.
\end{proof}

The next step in the construction is to show that in each of the nested subsets of $V(\xi_0,\eta_\ell)$ given by the previous lemma we will have not only the bases  $\alpha_{t^1,\ell}$ and $\alpha_{t^2,\ell}$, but also other components of $\gamma(t)$.  We split this step into two lemmas. The first lemma shows that for suitable choices of $A$-sequences the extended tails cut twice across the line Re$(z)=-\mu$ for $\mu>0$ arbitrarily large.

A finite block of $0$'s (respectively $0_{1}$'s or $0_{2}$'s) of length $k$ will be denoted by $0^k$ (respectively $0^k_1$ or $0^k_{2}$).

\begin{Lemma}\label{lemma:extended_tails}
Let $t$ be any $B$-sequence, and let $t^1, t^2\in \Sigma_A$ so that $\pi(t^1)=\pi(t^2)=t$. Assume also that $0_{1}t^1$ and $0_2t^2$ are allowable. Given any $\mu>0$, there exists $K>0$ such that for all $k>K$, there exist two continuous curves, denoted by $\tilde{\omega}_{1_20_1^k t^1}$ and
$\tilde{\omega}_{1_10_2^k t^2}$, that extend to infinity to the right and satisfy:
\begin{enumerate}
 \item[(a)] $\tilde{\omega}_{1_20_1^kt^1}$ and $\tilde{\omega}_{1_10_2^kt^2}$ are enlargements  (to the left) of the tails with  itineraries
 $1_20_1^kt^1$ and $1_10_2^kt^2$, respectively; and
 \item[(b)] $\tilde{\omega}_{1_20_1^kt^1}$ and $\tilde{\omega}_{1_10_2^kt^2}$ cuts twice across Re$\,z=-\mu$.
\end{enumerate}
\end{Lemma}

\begin{proof}
Consider any $A$-sequence  of the form $0_{1}^mt^1$ with $m>0$ large enough so the unique tail $\omega_{0_1^mt^1}$ is parametrized by $z=(x,h(x))$, $x > R$ (see Proposition \ref{prop:ConjugTails}). Lemma \ref{lem:ends-in-zeros} and its proof imply that $\omega_{0_1^mt^1}$ is  $\varepsilon$-close to $z=(x,0), \ x > R$ (that is, $|h(x)|<\varepsilon$ for all $x>R$). 

By pulling back $\omega_{0_1^mt^1}$ using $g_{a}^{0}$, we obtain a new curve which is precisely the extended tail $\tilde{\omega}_{0_1^{m+1}t^1}$. Moreover, since the positive real line is repelling, we observe that $\tilde{\omega}_{0_1^{m+1}t^1}$ is also $\varepsilon$-close to $z=(x,0)$ with $x > g_{a}^0(R)$. Successive pullbacks via $g_{a}^{0}$ allow us to find a first positive integer $r$ such that the extended tail $\tilde{\omega}_{0_1^{m+r}t^1}$ will be $\varepsilon$-close to  $z=(x,0), \ x > p_a$ where $p_a$ is the real fixed point in the boundary of $\mathcal A^*(-a)$.

At this stage we pullback again $\tilde{\omega}_{0_1^{m+r}t^1}$ using  $g_{a}^{1}$. The obtained curve must coincide with the extended tail $\tilde{\omega}_{1_{2}0_{1}^{m+r}t^1}$. By construction it is a curve extending to infinity to the right, extends far into the left half plane (since $\tilde{\omega}_{0_1^{m+r}t^1}$ is close to $z=0$) and lies close to $q_a$ (the preimage of $p_a$ in $\partial  \mathcal A^*(-a)$)  since $\tilde{\omega}_{0_1^{m+r}t^1}$ gets arbitrarily close to $z=p_a$. For a given $\mu$, if we choose $m$ large enough, this pullback construction guarantees that $\tilde{\omega}_{1_{2}0_{1}^{m+r}t^1}$ cuts across Re$(z)=-\mu$ twice, as desired.

Analogously, if we start the construction with $0_{2}^m t^2, \ m>0$, we get a similar result for the extended tail $\tilde{\omega}_{1_{1}0_{2}^{m+r}t^2}$.
\end{proof}

\begin{Proposition} \label{prop:passes}
Fix $\ell \geq 0$. Let $t$ be any $B$-sequence with $t^1$ and $t^2$ its associated $A$-sequences. Let $\tau$ be any finite block of binary symbols of length $n$ and denote by $s=\tau 110^kt$, and by $s^1,s^2$ its associated $A$-sequences. Then there exists $K>0$ such that for all $k>K$, the following statements hold.
\begin{enumerate}
\item[(a)] The forward image of $W_{n,\ell}^{1_{1}}$ cuts three times (twice far to the left and once far to the right) across the extended tail $\tilde{\omega}_{1_20_1^kt^1}$. In other words, the hair $\gamma(1_{1}1_20_1^kt^1)$ cuts three times across $W_{n,\ell}^{1_{1}}$.

\item[(b)] The hair $\gamma(s^1)$ cuts three times across
\begin{equation*}
g_{a}^{\tau_{0}}\circ\cdots\circ g_{a}^{\tau_{n-1}}\left(W_{n,\ell}^{1_{1}}\right)\ .
\end{equation*}
Moreover one of the connected components of  $\gamma(s^1)$ in 
the above expression contains $\alpha_{s^1,l}$.
\end{enumerate}
\noindent
Analogous results corresponding for $s^2$ also apply.
\end{Proposition}

\begin{proof} The set $f_{a}\left(W_{n,\ell}^{1_{1}}\right)$ is a large semi-annulus in the upper half plane and intersects $T_0 \cup T_1$ far to the right and far to the left. From Lemma \ref{lemma:extended_tails} we can choose $K$ such that for all $k>K$, the extended tail $\tilde{\omega}_{1_{2}0_{1}^{k}t^1}$ cuts across the semi-annulus $f_{a}\left(W_{n,\ell}^{1_{1}}\right)$ three times: twice far to the left and one far to the right. Consequently there must be three components of $\gamma(1_{1}1_20_1^kt^1)$ in $W_{n,\ell}^{1_{1}}$ that map into three components of the extended tail $\tilde{\omega}_{1_{2}0_{1}^{k}t^1}$ in $f_{a}\left(W_{n,\ell}^{1_{1}}\right)$.

By taking suitable pullbacks of $W_{n,\ell}^{1_{1}}$ that follow the string of symbols $\tau=(\tau_0,\ldots,\tau_{n-1})$ and applying Lemma  \ref{lemma:two-nested}, we can now conclude statement (b) of the present proposition.
\end{proof}

\begin{Theorem}\label{theorem:indecom}
Let $\tau$ be a finite block of binary symbols of length $n$. There exists an increasing sequence of integers $k_j,\ j\geq 1$ so for the $B$-sequence 
\begin{equation}
\T = \tau110^{k_1}110^{k_2}110^{k_3}\ldots, \label{eq:IndSeq}
\end{equation}
its associated $A$-sequences $\T^1$ and $\T^2$ determine
two distinct hairs $\gamma(\T^1)$ and $\gamma(\T^2)$ that limit upon themselves, thus becoming non-landing hairs.
\end{Theorem}

\begin{proof} 
Let $\ell>0$, $p_0=n+1$ and for each $i\geq 1$, set
 \begin{equation*}
p_i = n+1 + \sum\limits_{j=1}^{i} k_{j} + \ 2i.
 \end{equation*}
The symbol at the $p_i^{\text{th}}$ position in (yet to be constructed) $\T$ is equal to $1$. In general, the symbol at position $p_i$ in the $A$-sequence $\T^j,\ j=1,2$ can be either $1_1$ or $1_2$, depending on the previous entry in $\T^j$. In what follows we assume that all symbols at the $p_i^{\rm th}$ position in $\T^1$ are equal to $1_{1}$ (hence, all symbols at position $p_{i}^{\rm th}$ in $\T^2$ are $1_{2}$). It will become clear from the proof that this assumption is without lost of generality. Moreover, we shall prove only that  $\gamma(\T^1)$  limits upon itself since the case of $\gamma(\T^2)$ follows in the same way.

Let $s$ be any $B$-sequence not ending in all $0$'s, and let $s^1$ and $s^2$ be its associated $A$-sequences. For each $j\geq 1$ we aim to construct inductively $B$-sequences of the form 
\begin{eqnarray}\label{eq:sigma-j}
u_j &=& \tau 110^{k_1}110^{k_2}\ldots110^{k_j}\sigma^{p_j}(s),
 \end{eqnarray}
so that, by carefully selecting longer blocks of zeros, the $u_j$ will converge to the $B$-sequence $\T$ with the desired properties. Note that for all $j$, $u_j$ will be a concatenation of the first $p_j^{\rm th}$ symbols in $\T=(\mathfrak{t}_0, \mathfrak{t}_1, \ldots)$ with $\sigma^{p_j}(s)$. In other words,
\begin{eqnarray}\label{eq:terms}
u_j &=& (\mathfrak{t}_0, \mathfrak{t}_1, \ldots, \mathfrak{t}_{p_j-1}, s_{p_j}, s_{p_j+1},\ldots).
\end{eqnarray}
As before, $u_{j}^1$ and $u_{j}^2$ will denote the $A$-sequences that project into $u_j$.

We first show how to define $u_1$. Applying Proposition \ref{prop:passes}(a), we can choose an integer $k_1>0$ so the forward image of $W_{p_0-1,\ell}^{1_{1}}$ is cut across by the extended tail $\tilde{\omega}_{1_20_1^{k_{1}}\sigma^{p_1}(s^1)}$ three times (twice far to the left and once far to the right). Thus the hair $\gamma(1_{1}1_20_1^{k_{1}}\sigma^{p_1}(s^1))$ cuts three times across $W_{p_0-1,\ell}^{1_{1}}$. For the given block $\tau$, denote by $\tau_1$ its corresponding $A$-block that makes $u_1^1=\tau_1 1_1 1_2 0_1^{k_1}\sigma^{p_1}(s^1)$ an allowable $A$-sequence. Thus Proposition \ref{prop:passes}(b) implies that
$$  
g_{a}^{\mathfrak{t}_{0}}\circ\cdots\circ g_{a}^{\mathfrak{t}_{p_0-2}}\left(W_{p_0-1,\ell}^{1_{1}}\right)
$$
is a subset of $V\left(\xi_{0},\eta_{\ell}\right)$ that contains $\alpha_{u_{1}^1,\ell}$ and (at least) two other components  of the hair $\gamma(u_{1}^1)$.

For $j>1$ we apply again Proposition~\ref{prop:passes}(a) to the sequence
$$\tilde{\tau} 110^{k_j}\sigma^{p_j}(s),$$
where $\tilde{\tau}=\tau110^{k_1}110^{k_2}\ldots 110^{k_{j-1}}$ and with an integer $k_{j}>0$ large enough so that the forward image of $W_{p_{j}-1,\ell}^{1_{1}}$ is cut across by the extended tail $\tilde{\omega}_{1_20_1^{k_{j}}\sigma^{p_j}(s^1)}$ three times (twice far to the left and once far to the right). Thus the hair $\gamma(1_{1}1_20_1^{k_{j}}\sigma^{p_j}(s^1))$ cuts three times across $W_{p_{j}-1,\ell}^{1_{1}}$. Pulling back $W_{p_{j}-1,\ell}^{1_{1}}$ through the finite block of binary symbols given by $\tilde{\tau}$, we get from Proposition \ref{prop:passes}(b) that 
$$  
g_{a}^{\mathfrak{t}_{0}}\circ\cdots\circ g_{a}^{\mathfrak{t}_{p_j-2}}\left(W_{p_{j}-1,\ell}^{1_{1}}\right)
$$
is a subset of $V\left(\xi_{0},\eta_{\ell}\right)$ that contains $\alpha_{u_{j}^1,\ell}$ and (at least) two other components of the hair $\gamma(u_{j}^1)$.

As $j$ tends to infinity, the sequence $u_{j}$ converges to the desired sequence $\T$. Indeed, from Proposition \ref{prop:passes}(b), the sets 
$$  
g_{a}^{\mathfrak{t}_{0}}\circ\cdots\circ g_{a}^{\mathfrak{t}_{p_j-2}}\left(W_{p_{j}-1,\ell}^{1_{1}}\right)
$$
form a nested sequence of subsets of $V(\xi_{0},\eta_{\ell})$ that contains $\alpha_{u_{j},\ell}$ and two further components of the hair $\gamma(u_{j})$. Since $\ell\geq 0$ was selected in an arbitrary manner, as $j\to \infty$ we obtain $V\left(\xi_{0},\infty \right)$ contains $\omega_{\T^1}$ and infinitely many distinct components of the hair $\gamma(\T^1)$ accumulating on it.
To see this, let $z\in \omega_{\T^1}$ be given and select $\ell_z>0$ large enough so $z$ lies in $\alpha_{\T^1,\ell_z}$. As before, this base lies in the target $V(\xi_0,\eta_{\ell_z})$. Since $\T$ has been already constructed, we may now apply Proposition~\ref{prop:passes} to $\gamma(\T^1)$ itself. There exists an integer $N>0$ sufficiently large so for all $n>N$, $\gamma(\sigma^{p_n+1}(\T^1))$ cuts across $W_{p_n-1,\ell_z}^{1_1}$ in (at least) three components, say $A^n_1, A^n_2$ and $A^n_3$, with one of them being the component of the tail $\omega_{\sigma^{p_n+1}(\T^1)}$.

As $n$ increases, the diameter of $g_a^{\mathfrak{t}_0}\circ\ldots \circ g_a^{\mathfrak{t}_{p_n}}(W_{p_n-1,\ell_z}^{1_1})$ decreases, nonetheless, Proposition~\ref{prop:passes} implies that for each $j=1,2,3$,
$$ g_a^{\mathfrak{t}_0}\circ\ldots \circ g_a^{\mathfrak{t}_{p_n}}(A_j^n)~~\text{cuts across}~~ g_a^{\mathfrak{t}_0}\circ\ldots \circ g_a^{\mathfrak{t}_{p_n}}(W_{p_n-1,\ell_z}^{1_1}).$$
Since $g_a^{\mathfrak{t}_0}\circ\ldots \circ g_a^{\mathfrak{t}_{p_n}}(A_j^n)=\alpha_{\T^1,\ell_z}$ for some $j$, the pullback images of the remaining $A_j^n$'s accumulate lenghtwise over $\alpha_{\T^1,\ell_z}$ (and thus along $z$) for each $n>N$.

To see that $\gamma(\T^1)$ accumulates on each of its points and not only on its tail portion, we may perform the same construction for the sequences
$$
110^{k_i}110^{k_{i+1}}\ldots,
$$
for $i\geq 1$.
Then we may pullback the corresponding hairs and their accumulations by the appropriate inverse branches of $f_{a}$ to show that $\gamma(\T^1)$  must accumulate on any point in the hair $\gamma(\T^1)$. 
\end{proof}

\section{Geometry of hairs}
\label{section:indecom}

This last section will focus on the geometry of sets $\overline{\gamma(t)}$ when $t$ is a $B$-sequence that is either periodic or a sequence as constructed in Theorem~\ref{theorem:indecom}. We show that in the former case, the closure of the hair has a landing point in $\partial \A^*(-a)$ while, in the latter case, $\overline{\gamma(t^i)}$ is an indecomposable continuum for each $i=1,2$.

\begin{Proposition}\label{prop:PerSeqLand}
Let $t$ be a periodic $B$-sequence and $t^1, t^2$ their periodic associated $A$-sequences. Then, the hairs $\gamma(t^1)$ and $\gamma(t^2)$ land at two (repelling) periodic points $p_1,p_2\in \partial\A^*(-a)$.
\end{Proposition}

\begin{proof}

Let $t=\overline{t_0\ldots t_{n-1}}$, so every block of $0$'s in $t$ has bounded length. We work only with $t^1$ since the other case follows similarly.

There exists a value $\delta>0$ and a closed ball $\overline{B_{\delta}(0)}$ so the orbit of $\overline{\gamma(t^1)}$ stays always outside $\overline{B_{\delta}(0)}$. Analogously, we can find a real value $m$ so the orbit of $\overline{\gamma(t^1)}$ does not intersect the half plane Re$(z)< m<0$.

Let $R>0$ as in Proposition~\ref{prop:ConjugTails}. Then, there exists a value $M\geq R$ for which the intersection of the tail $\omega_{t^1}$ with the line Re$(z)=M$ is a single point. Since the orbit of the origin escapes along the positive real line, there exists an integer $N=N(M)>0$ for which $f_a^{N-1}(0)<M\leq f_a^{N}(0)$. Select $0<\epsilon\leq \delta$ small enough so for $B_\epsilon=\overline{B_\varepsilon(0)}$ and $f_a(B_\epsilon),\ldots, f^{N}_a(B_\epsilon)$, they are compact domains contained in $T_0$. Finally, select $m'\leq m<0$ so the left half plane Re$(z)\leq m'$ maps completely inside $B_\epsilon$.

With these constants we can define an open, connected, simply connected region $E$ as in the proof of Theorem~\ref{thm:BdOrbit}. Then, the pullback process that defines the hair $\gamma(t^1)$ together with the contracting map $\Psi_\ell= \psi_{t_0^1}^{-1}\circ(g_a^{t_0^1}\circ \cdots \circ g_a^{t_{n-1}^1})^\ell \circ \psi_{t_0^1}$ (where $(g_a^{t_0^1}\circ \cdots \circ g_a^{t_{n-1}^1})^\ell$ denotes the $\ell$-fold composition with itself), shows that $\overline{\gamma(t^1)}\setminus \gamma(t^1)$ is a unique periodic point in $E$ with itinerary $t^1$. By Theorem~\ref{thm:BdOrbit}, this periodic point lies in $\partial A^*(-a)$.
\end{proof}

\begin{Theorem}\label{thm:DoesNotSeparate}
Let $\T$ be a $B$-sequence as in Theorem~\ref{theorem:indecom}, $\T^1, \T^2$ be its associated $A$-sequences that project onto $\T$ and $\gamma(\T^1),\gamma(\T^2)$ their associated hairs. Then, the closure of each of these hairs is an indecomposable continuum.
\end{Theorem}

\begin{proof}
As before, we restrict the proof to the $A$-sequence $\T^1=\tau_1 1_1 1_2 0_1^{k_1}1_1 1_2 0_1^{k_2}\ldots$. Let $\Gamma^1$ be the closure of the hair $\gamma(\T^1)$. In order to apply Theorem~\ref{theo:curry}, we must verify first that $\Gamma^1$ does not separate the plane, as from Theorem~\ref{theorem:indecom} $\gamma(\T^1)$ is a curve that accumulates upon itself.

First observe that $\Gamma^1$ is a set with bounded negative real part. Indeed, the first block of $0_i$'s in $\T^1$ has finite length and this implies $\gamma(\T^1)$ (and thus $\Gamma^1$) lies to the right of the line Re$(z)=m$ for some $m<0$. For simplicity, denote by $k_0\in\{0_1,0_2,1_1,1_2\}$ the first entry in $\tau_1$, so $\Gamma^1$ lies inside the set $T_{k_0}$ and to the right of the line Re$(z)=m$. Hence
$$\bbC\setminus (T_{k_0}\cap \{z~|~\text{Re}(z)>m\})\subset \bbC\setminus \Gamma^1,$$
and since the bounderies of $T_0$ and $T_1$ are the graphs of strictly monotonic functions, the set in the left hand side is in fact a single component with unbounded imaginary part. This implies the existence of a single complementary component $U\subset \bbC\setminus \Gamma^1$ with unbounded imaginary part, while all other complementary components must have bounded imaginary parts. Also, note that $-a\in U$ and since $-a$ is a fixed point then $-a\in f_a^k(U)$ for all $k>0$.

Let $V \neq \emptyset $ be a component in $\bbC\setminus \Gamma^1$ with bounded imaginary part. Firstly, we assume that  $V\cap \J(f_a)\neq \emptyset$. Then,  by Montel's Theorem,  there exists $N>0$ such that $-a\notin f^k(V)$ for $0\leq k<N$ while $-a\in f^N_a(V)$, since $-a$ is not an exceptional value (indeed, $-a$ has infinitely many preimages).  But this implies that $f_a^{N-1}(V)$ is a component with bounded imaginary part inside $T_0\cup T_1$ that must contain a point in $f_a^{-1}(-a)$,  a contradiction since in $T_0 \cup T_1$ there are no preimages of $-a$ different from $-a$ itself. Secondly, we assume that  $V$ does not contain points in $\J(f_a)$ and thus is a Fatou component. Since the Fatou set coincides with the basin of attraction of $-a$, there exists an integer $N>0$ so that $f^N_a(U)=\A^*(-a)$ and thus $\partial \A^*(-a)\subset \partial U$, which is a contradiction with Proposition~\ref{prop:consequences_pol_like}(b), as only one point in $\partial \A^*(-a)$ has itinerary $\T^1$.

We conclude that $\Gamma^1$ does not separate the plane and by Theorem~\ref{theo:curry}, it is either an indecomposable continuum or the union of two indecomposable continua. Nevertheless, $\gamma(\T^1)$ has a unique tail extending to infinity so it cannot be the union of two indecomposable continua.
\end{proof}

\begin{Remark}
The above two results describe an interesting relationship between the combinatorics of an itinerary $s$ and the landing properties of the curve $\gamma(s)$. Indeed, if $s$ is an $A$-sequence whose blocks of $0_i$'s have bounded length, the same arguments as in Proposition~\ref{prop:PerSeqLand} can be applied to show that the accumulation set of the hair is no other than the unique point $p(s)\in \partial \A^*(-a)$ that follows the given itinerary. The key step is to ensure that pullbacks of the tail are always bounded away from the postcritical and asymptotical orbits (see, for instance, Proposition 3.6 in \cite{Fag}), and this is always the case for $A$-sequences whose blocks of $0_i$'s have bounded length.

In contrast, whenever $s=\T$ is a $B$-sequence as in Theorem~\ref{theorem:indecom}, then the corresponding point $p(\T^1)$ in $\partial \A^*(-a)$ is an accumulation point of $\gamma(\T^1)$, nevertheless, is not its endpoint.
\end{Remark}

\begin{Theorem}\label{theorem:relation_inde_basin}
Consider $\T$ and $\gamma(\T^1)$ as in Theorem~\ref{theorem:indecom}, set $\Gamma^1= \overline{\gamma(\T^1)}$. Then, there exists a unique point $p_1\in \partial \A^*(-a)$ such that $p_1\in \Gamma^1$ and $p_1$ has itinerary $\T^1$.
\end{Theorem}

\begin{proof}
Recall that $\T=\tau110^{k_1}110^{k_2}\ldots$. For each $n\geq1$ define two preperiodic sequences given by
\begin{eqnarray*}
s_n&=&\tau\overline{110^{k_1}\ldots 110^{k_n}0}, \\
r_n&=&\tau\overline{110^{k_1}\ldots 110^{k_n}1}.
\end{eqnarray*}
Clearly, their associated $A$-sequences satisfy $s^1_n < \T^1 < r^1_n$ for all $n\geq 1$ with respect to the distance induced by the usual order in $\Sigma_B$.
Moreover, if $p(s^1_n)$ and $q(r^1_n)$ denote endpoints of the corresponding preperiodic hairs associated to $s^1_n$ and $r^1_n$, then it follows by Proposition~\ref{prop:PerSeqLand} that $p(s^1_n), q(r^1_n)\in \partial \A^*(-a)$.
Let $p(\T^1)$ be the unique point in $\partial \A^*(-a)$ following the itinerary $\T^1$ under the action of $f_a|\partial\A^*(-a)$.\\

In the Euclidean distance restricted to $\partial \A^*(-a)$ we obtain $|p(s^1_n)-q(r^1_n)|\to 0$ as $n\to +\infty$ and clearly $p(s^1_n)<p(\T^1)<q(r^1_n)$ with the order inherited by $S^1$ under the continuous extension of the B\"ottcher mapping $\Phi_a:\bbD\to \A^*(-a)$.

For each $n$, $\gamma(\T^1)$ belongs to the region bounded above and below by $\gamma(r^1_n)$ and $\gamma(s^1_n)$, and the arc $[p(s^1_n),q(r^1_n)]\subset\partial \A^*(-a)$ containing $p(\T^1)$. Thus, if $\gamma(\T^1)$ accumulates on $\partial \A^*(-a)$, then it must accumulate at the point $p(\T^1)$.

In order to show that $\overline{\gamma(\T^1)}$ accumulates on the boundary of $\A^*(-a)$, we employ the polynomial-like construction and the symbolics of $\T^1$. Let $m>0$ and consider the hair associated to the itinerary
$$t_m=0_{k_m}11 0_{k_{m+1}}11 \ldots.$$

Note that $t_m$ is the image of $\T^1$ under some iterates of the shift $\sigma|{\Sigma_A}$. By Proposition~\ref{prop:passes}, the hair $\gamma(t_{m})$ is clearly close to the origin and also the hair $\gamma(1t_m)$ intersects a left half plane. Hence, a portion of $\gamma(1t_m)$ must intersect $V_a$ and its preimages under $f_a|V_a$ following $\tau110^{k_1}\ldots 0^{k_{m-1}}1$ must accumulate in $\partial \A^*(-a)$ since by Proposition~\ref{proposition:pol_like},
$$\bigcap_{n\geq 0} f_a^{-n}(V_a)=\partial\A^*(-a),$$
and $m>0$ has been taken arbitrarily large.
\end{proof}

We end this section by briefly discussing the set of points in the Julia set that follow non-binary sequences. Note first that Proposition~\ref{prop:ConjugTails} can be easily modified to show the existence of tails with itineraries corresponding to the partition $\cup_{j\in \bbZ} T_j$. If $s=(s_0,s_1, \ldots)$ is a sequence that contains infinitely many non-binary symbols and $T_{s_j}\subset f_a(T_{s_{j+1}})$ for all $j\geq 0$. Let $\omega_s$ be the tail associated to $s$ that lies in the right half plane Re$(z)>R$. If in addition $s$ has (if any) blocks of $0$'s with bounded length, then the pullbacks of $\omega_s$ are always bounded away from the postcritical and asymptotic orbits, hence $\gamma(s)$ is a landing hair with an endpoint $p(s)$. By Proposition~\ref{prop:consequences_pol_like}, $p(s)$ (and in fact $\overline{ \gamma(s)}$) do not lie in the boundary of any Fatou component, as otherwise, $s$ will have to end in a binary sequence. This establishes

\begin{Corollary}
If $s=(s_0, s_1,\ldots)\in \bbZ^\bbN$ is realizable by $f_a$, contains infinitely many non-binary symbols and its blocks of $0$'s have bounded lenght, then $\gamma(s)$ is a landing hair and a buried component of $\J(f_a)$.
\end{Corollary}

\subsection*{Acknowledgments}
The first and second author are both partially supported by the European network 035651-2-CODY, by MEC and CIRIT through the grants MTM2008Ð01486 and 2009SGR-792, respectively. The second author is also partially supported by MEC  through the grant MTM2006-05849/Consolider (including a FEDER contribution). The third  author is supported by CONACyT grant 59183, CB-2006-01. She would also like to express her gratitude to Universitat de Barcelona and Universitat Rovira i Virgili for their hospitality in the final stages of this article.

\bigskip
\noindent
{\small Antonio Garijo \& Xavier Jarque}\\
{\small Dept. d'Enginyeria Inform\`atica i Matem\`atiques}\\
{\small Universitat Rovira i Virgili}\\
{\small Av. Pa\"isos Catalans 26}\\
{\small Tarragona 43007, Spain}\\

\noindent
{\small M\'{o}nica Moreno Rocha}\\
{\small Centro de Investigaci\'on en Matem\'aticas}\\
{\small Callej\'on Jalisco s/n}\\
{\small Guanajuato 36240, Mexico}


\begin{thebibliography}{Abc}

\bibitem{Ba} K.~Bara\'nski, \emph{Trees and hairs for some hyperbolic
entire maps of finite order}. Math. Z. 257,  33-59 (2007).
\medskip
\bibitem{BJR} K.~Bara\'nski, X.~Jarque and L.~Rempe \emph{Brushing the hairs of transcendental entire functions}. Preprint
arXiv:1101.4209.
\bibitem{Berg} W.~Bergweiler, \emph{Iteration of meromorphic functions.}  
Bull. Amer. Math. Soc. (N.S.) 29, no. 2, 151Ð188 (1993).
\medskip
\bibitem{C} S.~Curry, \emph{One-dimensional nonseparating plane continua with
disjoint $\epsilon$-dense subcontinua}.  Topol.~and its Appl. 39, 145-151 (1991).
\medskip
\bibitem{D}  R.~L.~Devaney, \emph{Knaster-like continua and complex dynamics}. Erg. Th. and Dyn. Sys. 13, 627-634 (1993).
\medskip
\bibitem{DK}  R.~L.~Devaney and M. Krych, \emph{Dynamics of {\rm exp}$(z)$}. Erg. Th. and Dyn. Sys. 4, 35-52 (1984).
\medskip
\bibitem{DJ1} R.~L.~Devaney and X.~Jarque, \emph{Misiurewicz points for
complex exponentials}. Int. J. Bifurc. and Chaos. 7, 1599-1616 (1997).
\medskip
\bibitem{DJ2} R.~L.~Devaney and X.~Jarque, \emph{Indecomposable continua in complex dynamics}. Conform. Geom. Dyn. (electronic), 1-12 (2002).
\medskip
\bibitem{DJM} R.~L.~Devaney, X.~Jarque and M.~Moreno Rocha, \emph{Indecomposable continua and Misiurewicz points in exponential dynamics}. Int. J. Bifurc. and Chaos. 15(10), 3281-3294  (2005).
\medskip
\bibitem{DT} R.~L.~Devaney and F. Tangerman, \emph{Dynamics of entire
functions near the essential singularity}. Ergod. Th. and Dynam. Sys.  6, 489-503  (1986).
\medskip
\bibitem{DH} A.~Douady and J.~H.~Hubbard,\emph{On the dynamics of
polynomial-like mappings}. Ann. Sci. \'Ecole Norm. Sup. 4, 
no. 2, 287-343  (1985).
\medskip
\bibitem{E} A.~\`E.~Er\"emenko, \emph{On the iteration of entire functions}. In
 Dynamical systems and ergodic theory (Warsaw, 1986), Banach Center
 Publ., Warsaw,  339-345 (1989).
\medskip
\bibitem{EL2} A.~\`E.~Er\"emenko and M.~Yu.~Lyubich, \emph{Dynamical properties of some classes of entire functions}. Ann.~Inst.~Fourier (Grenoble) 42 no.~4, 989-1020 (1992).
\medskip
\bibitem{Fag} N.~Fagella, \emph{Dynamics of the complex standard family}.
J. Math. Anal. Appl. 229 no.~1, 1-31 (1999). 
\medskip
\bibitem{F} P.~Fatou, \emph{Sur l'It\'eration des fonctions transcendentes
enti\`eres}. Acta Math. 47, 337-370  (1926).
\medskip
\bibitem{GK} L.~Goldberg and L.~Keen, \emph{A finitness theorem for a
dynamical class of entire functions}. Erg. Th. and Dyn. Syst.  6, 183-192  (1986).
\medskip
\bibitem{K} C.~Kuratowski, \emph{Topologie}. Volume {II}. Pa\'nstwowe Wydawnictwo Naukowe, Warzawa (1961). 
\medskip
\bibitem{MR} M.~Moreno Rocha, \emph{Existence of indecomposable continua for
unstable exponentials}. Top. Proc. 27, 233-244  (2002).
\medskip
\bibitem{Moro} S.~Morosawa, \emph{Local connectedness of Julia sets for transcendental entire functions}. In Proceedings of the International Conference on Nonlinear Analysis and Convex Analysis,  pp. 266-273. World Scientific, (1999).
\medskip
\bibitem{MU} S.~Morosawa, Y.~Nishimura, M.~Taniguchi, and T.~Ueda, \emph{Holomorphic dynamics}. Cambridge University Press, Cambridge (2000).
\medskip
\bibitem{Na} S.~B.~Nadler Jr., \emph{Continuum theory}. Marcel Dekker Inc., New York (1992).
\medskip
\bibitem{R} L.~Rempe, \emph{Rigidity of escaping dynamics for transcendental entire maps}. Acta Math. 203, 235-267 (2009). 
\medskip
\bibitem{R1} L.~Rempe, \emph{On nonlanding dynamic rays of exponential maps}. Ann. Acad. Sci. Fenn. Math. 32, 353-369  (2007).   
\medskip
\bibitem{R2} L.~Rempe, \emph{Siegel disks and periodic rays of entire
functions}. J. Reine Angew. Math. 624, 81-102 (2008).
\medskip
\bibitem{SZ} D.~Schleicher and J.~Zimmer,  \emph{Escaping Points of Exponential Maps}. J. London Math. Soc.  67 (2), 380-400 (2003).
\medskip
\bibitem{Why} G.~T.~Whyburn, \emph{Topological characterization of the Sierpi\'nski curve}. Fund. Math. 45, 320-324 (1958).

\end{thebibliography}
\end{document}